%% file: manuscript.tex
\documentclass[final,3p,times]{elsarticle}

\graphicspath{{figures/}}
\usepackage{tikz}
\usetikzlibrary{calc, patterns, decorations.pathmorphing, decorations.markings}

\usepackage{hyperref}
\usepackage{amsmath}
\usepackage{cleveref}

\usepackage{array} 
\usepackage{booktabs} 
\usepackage{multicol} 
\usepackage{algorithm, algorithmicx, algpseudocode}
\usepackage{xspace}
\usepackage{bm}

\usepackage{amsthm}
\usepackage{physics}
\usepackage{subcaption}
\usepackage{color}
\usepackage{import}
\graphicspath{{figures/}}
\usepackage{booktabs} 

\makeatletter
\def\ps@pprintTitle{}
\makeatother

\usepackage{pgfplots}
\usetikzlibrary{decorations.markings} 
\usepackage{pgfplotstable}
\pgfplotsset{compat=1.18}
\usepgfplotslibrary{groupplots}

\biboptions{authoryear}

\begin{document}

\begin{frontmatter}

\title{First-order equivalent static loads for dynamic response structural optimization}


\author[mymainaddress1]{Mordechay Buzaglo}
\ead{mordechayb@campus.technion.ac.il}

\author[mymainaddress2]{Nicol\`{o} Pollini}
\ead{nicolo@technion.ac.il}

\address[mymainaddress1]{Faculty of Mechanical Engineering, Technion - Israel Institute of Technology, Haifa, Israel}
\address[mymainaddress2]{Faculty of Civil and Environmental Engineering, Technion - Israel Institute of Technology, Haifa, Israel}

\begin{abstract}
A novel first-order equivalent static loads approach for optimization of structural dynamic response, F-ESL, is presented and compared to the basic equivalent static load formulation, ESL. F-ESL simplifies dynamic optimization problems by converting them into a series of static optimization sub-problems. The ESL algorithm in its original formulation does not have a guaranteed capability of reaching, or recognizing, final designs that satisfy necessary first-order optimality conditions.  
F-ESL addresses this limitation by including first-order terms directly into the equivalent static load definition. This new mathematical information guides the optimization algorithm more effectively toward solutions that satisfy both feasibility and optimality conditions. 
Using reproducible numerical examples, we show that F-ESL overcomes the known limitations of the original ESL, often with few outer function evaluations and fast convergence. At the same time, F-ESL maintains ESL simplicity, robustness, and ease of implementation, providing practitioners with an effective tool for structural dynamic optimization problems. 

\end{abstract}


\end{frontmatter}



\section{Introduction}
Optimization of the dynamic response plays a crucial role in the design of several engineering structures. 
For example, airplanes \citep{kennedy2014parallel}, cars \citep{torstenfelt2006structural}, bridges \citep{zordan2010bridge,martins2019optimization}, wind turbines \citep{pavese2017aeroelastic,blasques2012multi,hermansen2024multi}, buildings \citep{marzok2022topology,idels2023optimization}, and, on a smaller scale, microdevices such as MEMS resonators \citep{comi2016non,giannini2020topology,giannini2022topology}. 
The approaches for optimization of structural dynamic response can be classified into three groups, namely: 1) optimization of free vibration response \citep{diaaz1992solutions,pedersen2000maximization,ferrari2018eigenvalue}, 2) optimization of forced vibration response \citep{olhoff2016generalized,niu2018objective,li2021reduced}, and 3) optimization of transient time-dependent response \citep{pedersen2004crashworthiness,james2015topology,behrou2017topology}. 
Here, we discuss the design optimization of structural systems subjected to transient time-dependent dynamic excitation. 
In the time domain, the computational cost of structural optimization with dynamic loading is quite significant. The reason is the repetition of several time-history and sensitivity analyses during the design updates of the optimization process. 
High computational cost is still a limiting factor for practical engineering applications. 
For this reason, in the literature, some authors have focused on developing optimization approaches and strategies to reduce the computational complexity and cost of structural optimization for transient dynamic problems. 

In dynamic structural optimization procedures, the response and sensitivity analyses must be repeated many times during the design optimization process due to changes in the design variables. The high computational cost involved in reanalysis is one of the main obstacles in the optimization of large-scale structures. This is even more true for time-dependent dynamic problems. 
To alleviate this difficulty, Kirsch and co-authors developed the combined approximation approach (CA) \citep{kirsch2010reanalysis}. 
CA consists of a uniﬁed approach for efficient reanalysis and sensitivity reanalysis of linear, nonlinear, static and dynamic systems. 
In essence, for optimization of linear static response the CA approach allows one to approximate the displacements of the structure associated with a current design update (reflected in an updated system stiffness matrix) based on a reduced vector basis. The vector basis is used to reduce the dimensionality of the system matrices, thus reducing the computational cost associated with the solution of the system of equilibrium equations. 
Initially developed for optimization of linear static response \citep{kirsch1991reduced}, the CA approach has been applied to a variety of different cases, such as structural optimization of nonlinear response \citep{amir2008efficient}, vibration problems \citep{kirsch2004procedures}, linear \citep{kirsch2007efficient} and nonlinear \citep{kirsch2006nonlinear} dynamic response. 

A different strategy for reducing the computational cost during optimization of structural dynamic response is the use of reduced-order models. 
For example, the use of a reduced-order model for topology optimization of transient linear dynamic response is discussed in \cite{long2021topology}. The second-order Arnoldi
reduction (SOAR) scheme is employed to reduce the dimensions of the system matrices, and hence the computational cost during transient dynamic analysis. 
The use of reduced-order models for topology optimization of linear structural dynamic response is also the focus of the work of Zhao and Wang \citep{zhao2016dynamic}. Two approaches are used to reduce the order of the system matrices, the mode displacement method (MDM) and the mode acceleration method (MAM). 
Very recently, an optimization approach that combines Kirsch's CA approach and equivalent static loads (ESL) has been proposed in \cite{li2023efficient}. 
The computational cost in the optimization process is reduced in two ways. 
Reanalysis techniques (in the specific, CA \citep{kirsch2010reanalysis}) are used to approximate the structural response during repeated transient linear dynamic analyses. 
The design updates are then obtained by means of an ESL algorithm, which transforms the initial topology optimization problem for transient dynamic response into a topology optimization problem for linear static response with multiple loading conditions, one for each time step.

It can be seen that in the literature there are approaches for reducing the computational cost during optimization of time-dependent structural dynamic response.  
In this regard, equivalent static load (ESL) optimization approaches, if properly designed, are an effective tool for reducing the computational cost of structural optimization procedures of time-dependent dynamic response. 
ESL approximates dynamic or non-linear optimization problems with a sequence of static linear optimization sub-problems. In the case of dynamic response optimization, each time step translates into a vector of equivalent static loads. The loads cause displacements in the static linear system equivalent to those observed in the corresponding time step of the original dynamic problem. These loads are treated as multiple load cases in the static response optimization sub-problem.
In case of linear time-dependent dynamic response, the ESL are the loads $\mathbf{f}_{\text{eq},k}$ that in each $k$-th time-step produce displacements in the linear static structural system (defined by $\mathbf{K}_{L}$) equal to those obtained with transient time-history analysis $\mathbf{u}_{k}$:
\begin{equation}
    \mathbf{f}_{\text{eq},k} = \mathbf{K}_{L} \, \mathbf{u}_{k},
\end{equation}
where the displacements $\mathbf{u}_{k}$ are computed solving numerically the following dynamic structural analysis problem:
\begin{equation}
\label{eq:motion}
\begin{split}
    & \mathbf{M}\,\ddot{\mathbf{u}}_k + \mathbf{C}\,\dot{\mathbf{u}}_k + \mathbf{K}\,\mathbf{u}_k = \mathbf{p}_k \text{ for } k=1,\dots,n_t \\
    & \mathbf{u}_0 = \mathbf{0}, \; \dot{\mathbf{u}}_0=\mathbf{0}.
\end{split}
\end{equation}
In Eq.~\eqref{eq:motion}, $\mathbf{M}$ is the structural mass matrix, $\mathbf{K}$ is the stiffness matrix and $\mathbf{C}$ is the damping matrix and are computed, for example, by means of the finite element method. The vector $\mathbf{p}_k$ contains the time-dependent loads. Similar considerations can be extended to the case of nonlinear transient dynamic response.
In the literature, ESL approaches are a well-known approach for dynamic-response structural optimization, first introduced by \cite{choi2002structural}. 
 
From a chronological point of view, an ESL algorithm is first introduced for optimization of linear time-dependent dynamic response based on modal analysis in \cite{choi1999transformation}. The ESL algorithm proposed in \cite{choi1999transformation} is further discussed in \cite{kang2001structural,choi2002structural}.  
ESL is also applied to nonlinear dynamic response optimization for the design of car roof \citep{jeong2008structural} and frontal \citep{jeong2010non} structural components. 
Moreover, in \cite{park2005structural} the ESL algorithm is applied to shape optimization of linear structures subjected to dynamic loads. 
The application of ESL to the optimization of multi-body dynamic systems is discussed in \cite{tromme2016equivalent,tromme2018system}.
Recently, topology optimization of structures with nonlinear dynamic behavior has been discussed in \cite{lee2015nonlinear}.
The optimization of a front hood for nonlinear dynamic structural response based on the ESL algorithm is discussed in \cite{karev2017comparison,karev2019free}. 
The ESL algorithm has been applied also to the car roof crash test optimization problem in \cite{choi2018structural}.
The ESL algorithm in its original formulation is a zeroth-order approximation of the original problem and generates equivalent static loads that are independent of the design variables. Hence, it does not have a guaranteed capability of reaching, or recognizing, final designs that satisfy the necessary first-order optimality conditions. 
This limitation has already been discussed in the literature for the case of linear transient structural optimization problems \citep{stolpe2014equivalent, stolpe2018equivalent}.  
The goal of this paper is to bridge this gap in the context of optimization of linear dynamic response optimization.
An alternative to the original formulation of ESL for nonlinear transient structural optimization problems, a difference-based extension of the ESL algorithm (DiESL), has been recently proposed in \cite{Triller2021}. 
The main idea in DiESL is to compute incremental equivalent static loads by considering incremental displacements and the linear analysis structural stiffness matrix evaluated at the displaced structural configuration in each time-step. DiESL is further extended in \cite{triller2022difference} by considering heterogeneous time steps, adaptively placed at points in time where nonlinearities are dominant. 
Nevertheless, also DiESL considers only zeroth-order definitions of the ESL.
It should be noted that for the case of linear dynamic structural optimization problems, DiESL is equivalent to the original ESL algorithm \citep{stolpe2023diesl}.   
Recently, Stolpe and Pollini \citep{stolpe2023first} presented preliminary results on first-order equivalent static loads for static response optimization of 2D geometrically nonlinear truss problems. Their approach successfully modifies the original ESL approach using a first-order approximation of the equivalent loads with respect to the design variables. 
In the following, we present a first study of F-ESL for optimization of structural dynamic response. The approach is first presented and discussed. To demonstrate the performance of the proposed F-ESL method, a sequence of reproducible numerical examples is presented and discussed in detail. We begin with a simple illustrative problem, a one-degree-of-freedom two-bar system, to validate the theoretical findings and compare with the classical ESL approach. We then apply the proposed F-ESL approach to two additional design optimization problems: a two-story shear frame subjected to seismic excitation with displacement constraints, and a truss structure subjected to transient loading with stress and buckling constraints. The numerical examples highlight the ability of F-ESL to identify optimal design solutions with modest computational resources, while retaining the algorithmic simplicity of the original ESL approach. 

In what follows, Section \ref{sec:probform} presents the original and approximate optimization problem formulations; Section \ref{sec:senanal} describes in detail the derivation for the displacement response sensitivity with respect to the design variables, which is used in Section \ref{subsec:kktcondcomp} to prove that the proposed F-ESL approach is indeed a first-order approximation of the original problem. Several reproducible numerical examples are discussed in detail in Section \ref{sec:numres}, and final conclusions are drawn in Section \ref{sec:end}.

\section{Original and approximate problem formulations}
\label{sec:probform}
In the following section, we present a general formulation of an optimization problem for structures subjected to dynamic loads. We then discuss the details of the basic ESL and the proposed F-ESL approaches. Note that the general formulation used here can be used to describe any structural optimization problem with displacement-dependent objective function and constraints.

\subsection{Dynamic response optimization problem}
\label{subsec:dynresoptprob}

The optimization problem considered here consists of an objective function $f_c(\mathbf{x}, \mathbf{u}_k(\mathbf{x}))$ and $m$ constraints $g_j(\mathbf{x}, \mathbf{u}_k(\mathbf{x}))$ with $j=1,2,\dots,m$, all potentially dependent on the vector of the design variables and the system displacements at each $k$-th time step. 
The problem is presented here for an arbitrary number of design variables $n$, time steps $n_t$, and degrees of freedom (DOFs). The reference dynamic response optimization problem formulation reads:
\begin{equation}
\label{eq:dynoptform}
\begin{split}
    \underset{\mathbf{x} \in {\rm I\!R}^n}{\text{minimize }} &  f_c(\mathbf{x},\mathbf{u}_k (\mathbf{x}))\\
    \text{subject to:  } & g_j(\mathbf{x},\mathbf{u}_k (\mathbf{x})) \leq 0, \hspace{5mm} j=1,2,\dots,m \\
    & x_i^{min} \leq x_i \leq x_i^{max}, \hspace{5mm} i=1,2,\dots,n \\
    \text{with:  } &\mathbf{M}(\mathbf{x}) \mathbf{\ddot{u}}_k(\mathbf{x})
    +\mathbf{C}(\mathbf{x}) \mathbf{\dot{u}}_k(\mathbf{x}) 
    +\mathbf{K}(\mathbf{x}) \mathbf{u}_k(\mathbf{x}) = \mathbf{f}_k\\
    & \text{for } k=1,2,\dots,n_t
\end{split}
\tag{$\mathcal{P}_D$}
\end{equation}
where $\mathbf{x}$ is the vector of the design variables of length $n$; 
$\mathbf{u}_k (\mathbf{x})$, $\dot{\mathbf{u}}_k (\mathbf{x})$ and $\ddot{\mathbf{u}}_k (\mathbf{x})$ represent the displacement, velocity and acceleration vectors of all free DOFs at time-step $k$; 
$k\in[1,n_t]$ is the full range of time steps over which the dynamic response is computed; $k\in[s_t,n_t] \subset [1,n_t]$, instead, is the discrete time window in which the objective function and constraints are evaluated. 
The design variables $x_i$ for $i=1,\dots,n$ are allowed to vary between $x_i^{min}$ and $x_i^{max}$. 
Moreover, $\mathbf{M} (\mathbf{x})$, $\mathbf{C} (\mathbf{x})$ and $\mathbf{K} (\mathbf{x})$ are the global mass, damping, and stiffness matrices of the system. The vector $\mathbf{f}_k$ represents the external force acting on each free DOF at time step $k$.

\subsection{Equivalent static loads formulations}
\label{subsec:eslmeth}
The ESL optimization algorithm replaces the original dynamic response optimization problem with a sequence of static response optimization sub-problems, which are solved iteratively until predefined stopping criteria are met. 
In every outer iteration, the dynamic response of the structure is evaluated and the equivalent static loads are computed. Then, inner iterations are performed to optimize the equivalent linear static system, subjected to multiple equivalent static load cases (one for each time step). 
The general formulation of the optimization sub-problems in the $W$-th outer iteration is defined in \eqref{eq:statoptform}. 
This formulation applies to both the ESL and the F-ESL approaches:
\begin{equation}
\label{eq:statoptform}
\begin{split}
    \underset{\mathbf{x} \in {\rm I\!R}^n}{\text{minimize }} &  f_c(\mathbf{x}, \tilde{\mathbf{u}}_k (\mathbf{x}))\\
    \text{subject to:  } & g_j(\mathbf{x}, \tilde{\mathbf{u}}_k (\mathbf{x})) \leq 0, \hspace{5mm} j=1,\dots,m \\
    & x_i^{min} \leq x_i \leq x_i^{max}, \hspace{5mm} i=1,\dots,n \\
    \text{with:  } & \mathbf{K} (\mathbf{x}) \tilde{\mathbf{u}}_k(\mathbf{x}, \mathbf{x}_W) 
     = \mathbf{f}^{ESL}_{k} (\mathbf{x}_W)\\ 
    \text{or } & \mathbf{K} (\mathbf{x}) \tilde{\mathbf{u}}_k(\mathbf{x}, \mathbf{x}_W) 
     = \mathbf{f}^{FESL}_{k} (\mathbf{x}, \mathbf{x}_W)\\ 
     &\text{for } k=s_t,\dots,n_t.
\end{split}
\tag{$\mathcal{P}_{S,W}$}
\end{equation}
In the optimization problem \eqref{eq:statoptform}, $\mathbf{x}$ is the optimization variable vector; $\tilde{\mathbf{u}}_k$ is the displacement vector of the static sub-problem; $\mathbf{f}^{ESL}_{k}$ and $\mathbf{f}^{FESL}_{k}$ are the ESL and F-ESL equivalent loads defined in Eq.~\eqref{eq:esldef} and ~\eqref{eq:feslfdef}. 
Moreover, $\mathbf{x}_W$ is the vector with the current values of the optimization variables obtained in \eqref{eq:statoptform} in the $W$-th iteration. 
It is used to perform the dynamic analysis and to recompute the corresponding equivalent static loads $\mathbf{f}^{ESL}_{k}$ or $\mathbf{f}^{FESL}_{k}$ that are then passed as starting point to the subsequent $W+1$-th sub-problem \eqref{eq:statoptform}. 

The basic ESL approach is based on a zeroth-order formulation. It means that with ESL the equivalent static loads do not depend on the design variables $\mathbf{x}$ in the inner static response optimization sub-problems. 
Therefore, with ESL the equivalent loads are constant in each sub-problem. The loads in the ESL formulation are defined as follows:
\begin{equation}
\label{eq:esldef}
    \mathbf{f}^{ESL}_{k}(\mathbf{x}_W) = \mathbf{f}^{eq}_{k}(\mathbf{x}_W) =\mathbf{K}(\mathbf{x}_W) \mathbf{u}_k(\mathbf{x}_W).
\end{equation}
Thus, for each time-step $k=1,\dots,n_t$ , the equivalent static loads used in \eqref{eq:statoptform} are $\mathbf{f}^{ESL}_{k} (\mathbf{x}_W)$.

In this paper, we suggest to rely on a novel formulation for the equivalent static loads, namely F-ESL. 
This formulation includes a first-order approximation of the equivalent static loads:
\begin{equation}
\label{eq:feslfdef}
\begin{split}
&\mathbf{f}^{FESL}_{k}(\mathbf{x},\mathbf{x}_W) =  
 \mathbf{f}^{eq}_{k}(\mathbf{x}_W) +\nabla \mathbf{f}^{eq}_{k}(\mathbf{x}_W) \left( \mathbf{x}-\mathbf{x}_W \right),
\end{split}
\end{equation}
where 
\begin{equation}
\label{eq:nablafESL}
\begin{split}
\nabla \mathbf{f}^{eq}_{k}(\mathbf{x}) &= \nabla \left( \mathbf{K}(\mathbf{x}) \mathbf{u}_k(\mathbf{x}) \right)\\
&= \nabla{\mathbf{K}(\mathbf{x})} \mathbf{u}_k(\mathbf{x}) + \mathbf{K}(\mathbf{x})\nabla{ \mathbf{u}_k(\mathbf{x})} .
\end{split}
\end{equation}
Thus
\begin{equation}
\label{eq:feslfdef2}
\begin{split}
&\mathbf{f}^{FESL}_{k}(\mathbf{x},\mathbf{x}_W) = 
\mathbf{K}(\mathbf{x}_W) \mathbf{u}_k(\mathbf{x}_W) + \\
& + \left[ \nabla{\mathbf{K}(\mathbf{x}_W)} \mathbf{u}_k(\mathbf{x}_W) + \mathbf{K}(\mathbf{x}_W)\nabla{ \mathbf{u}_k(\mathbf{x}_W)} \right] \,\left( \mathbf{x}-\mathbf{x}_W \right). 
\end{split}
\end{equation}
In the proposed F-ESL approach, for each time-step $k=s_t,\dots,n_t$, the equivalent static loads considered in \eqref{eq:statoptform} are $\mathbf{f}^{FESL}_{k} (\mathbf{x},\mathbf{x}_W)$. 

With both the ESL and F-ESL approaches, after solving each static response optimization sub-problem, the updated design variable vector $\mathbf{x}_W$ is used to perform the dynamic analysis, compute new zeroth- or first-order equivalent static loads, and proceed to the subsequent static response optimization sub-problem using  $\mathbf{x}_W$ as the starting point. 
This procedure is repeated until the following stopping criterion is satisfied:
\begin{equation}
\label{eq:stoppingcond}
    ||\mathbf{x}_{W+1}-\mathbf{x}_{W}|| \leq \epsilon
\end{equation}
where $\epsilon$ is a small positive number (e.g., $10^{-6}$). The pseudo codes for both F-ESL and ESL formulations are presented in Algorithm~\ref{pseudoalgo}.

\begin{algorithm}
    \centering
    \include{elements/pseudo_algo}
    \caption{\small Pseudo code for the Equivalent Static Load (ESL) and First-order ESL (F-ESL) formulations. The algorithm iteratively updates the design variables by performing dynamic-response analysis, computing the equivalent static loads, solving the linear static optimization sub-problem, until stopping criteria are satisfied.}
    \label{pseudoalgo}
\end{algorithm}

\section{Sensitivity analysis}
\label{sec:senanal}

In general, both the objective function and the constraints depend on the system response as well as explicitly on the design variables. The displacement sensitivities, $\frac{\partial \mathbf{u}_k}{\partial x_i}$, are governed by the governing equations of the problem and, consequently, by the specific formulation of the equivalent static loads employed. Therefore, to fully define the gradients of the objective and constraint functions, the displacement derivatives must be evaluated separately for each formulation considered.

The equivalent static loads are evaluated using the dynamic response of the original problem, after performing the time-history analysis. 
To this end, first the response in each time-step is evaluated using a numerical integration scheme (Newmark-$\beta$ \citep{chopra2007dynamics}, for example). The response then can be expressed as: 
\begin{equation}
\label{eq:u_dyn}
\begin{split}
\mathbf{u}_k(\mathbf{x}) = &\mathbf{K}^{-1}(\mathbf{x}) \left( -\mathbf{M}(\mathbf{x}) \ddot{\mathbf{u}}_k(\mathbf{x})  -\mathbf{C}(\mathbf{x}) \dot{\mathbf{u}}_k(\mathbf{x})+
\mathbf{f}_k \right) \\
=&\mathbf{K}^{-1}(\mathbf{x}) \mathbf{f}^{eq}_{k}(\mathbf{x}).
\end{split}
\end{equation}
The derivative of the displacements is computed as follows:
\begin{equation}
\label{eq:du_dx}
\begin{split}
\frac{\partial \mathbf{u}_k(\mathbf{x})}{\partial x_i} = \frac{\partial \mathbf{K}^{-1}(\mathbf{x})}{\partial x_i} \mathbf{f}^{eq}_{k}(\mathbf{x}) + \mathbf{K}^{-1}(\mathbf{x}) \frac{\partial \mathbf{f}^{eq}_{k}(\mathbf{x})}{\partial x_i}
\end{split}
\end{equation}
Considering Eq.~\eqref{eq:u_dyn}, the derivative of the equivalent static loads can be expressed as follows:
\begin{equation}
\label{eq:dfESL_dx}
\begin{split}
\frac{\partial \mathbf{f}^{eq}_{k}(\mathbf{x})}{\partial x_i} = &-\frac{\partial \mathbf{M}(\mathbf{x})}{\partial x_i} \ddot{\mathbf{u}}_k(\mathbf{x}) - \mathbf{M}(\mathbf{x}) \frac{\partial \ddot{\mathbf{u}}_k(\mathbf{x})}{\partial x_i} \\
&- \frac{\partial \mathbf{C}(\mathbf{x})}{\partial x_i} \dot{\mathbf{u}}_k(\mathbf{x}) - \mathbf{C}(\mathbf{x}) \frac{\partial \dot{\mathbf{u}}_k(\mathbf{x})}{\partial x_i}.
\end{split}
\end{equation}
In contrast to what is stated in \cite{park2003validation}, in general, $\nabla \mathbf{f}^{eq}_{k}(\mathbf{x})$ (defined in Eq.~\eqref{eq:nablafESL}) is not necessarily zero \textit{at an optimal solution $\mathbf{x}=\mathbf{x}^*$ of the original problem}, as shown also numerically in the first example of Sec.~\ref{sec:numres}. 
In the basic ESL approach, the $\mathbf{f}^{ESL}_{k}(\mathbf{x}_W)$ vector is constant during the static-response optimization sub-problem \eqref{eq:statoptform}, since it does not depend directly on $\mathbf{x}$. Therefore, 
its derivative with respect to a design variable $x_i$ using Eq.~\eqref{eq:esldef} is equal to zero: 
\begin{equation}
\label{eq:dfesl_dxtilde}
\begin{split}
\frac{\partial \mathbf{f}^{ESL}_k(\mathbf{x}_W)}{\partial x_i}=\frac{\partial \mathbf{f}^{eq}_k(\mathbf{x}_W)}{\partial x_i}=\mathbf{0}.
\end{split}
\end{equation}
The displacement gradient can then be calculated as follows:
\begin{equation}
\label{eq:duESL_tilde_dx}
\begin{split}
\frac{\partial \tilde{\mathbf{u}}^{ESL}_k(\mathbf{x})}{\partial x_i} &= \frac{\partial \mathbf{K}^{-1}(\mathbf{x})}{\partial x_i} \mathbf{f}^{ESL}_{k}(\mathbf{x}_W) + \mathbf{K}^{-1}(\mathbf{x}) \underbrace{\frac{\partial \mathbf{f}^{ESL}_{k}(\mathbf{x}_W)}{\partial x_i}}_{=0} \\
& =\frac{\partial \mathbf{K}^{-1}(\mathbf{x})}{\partial x_i} \mathbf{f}^{ESL}_{k}(\mathbf{x}_W).
\end{split}
\end{equation}

With the proposed F-ESL definition, where
$\mathbf{f}^{FESL}_{k}$ depends on $\mathbf{x}$, 
the displacement sensitivities are computed as follows:
\begin{equation}
\label{eq:duFESLf_tilde_dx}
\begin{split}
&\frac{\partial \tilde{\mathbf{u}}^{FESL}_k(\mathbf{x}, \mathbf{x}_W)}{\partial x_i} = \\
&\frac{\partial \mathbf{K}^{-1}(\mathbf{x})}{\partial x_i} \mathbf{f}^{FESL}_{k}(\mathbf{x},\mathbf{x}_W)
+ \mathbf{K}^{-1}(\mathbf{x}) \frac{\partial \mathbf{f}^{FESL}_{k}(\mathbf{x},\mathbf{x}_W)}{\partial x_i}  \\  &=\frac{\partial \mathbf{K}^{-1}(\mathbf{x})}{\partial x_i} \left( \mathbf{f}^{eq}_{k}(\mathbf{x}_W)+  \nabla{\mathbf{f}^{eq}_{k}(\mathbf{x}_W)} \left( \mathbf{x}-\mathbf{x}_W \right) \right) \\ &+ \mathbf{K}^{-1}(\mathbf{x}) \frac{\partial \mathbf{f}_k^{eq}(\mathbf{x}_W)}{\partial x_i}.
\end{split}
\end{equation}
Rearranging Eq.~\eqref{eq:duFESLf_tilde_dx} while considering Eq.~\eqref{eq:nablafESL} and Eq.~\eqref{eq:esldef}, we get the following expression:
\begin{equation}
\label{eq:duFESLf_tilde_dx_2}
\begin{split}
&\frac{\partial \tilde{\mathbf{u}}^{FESL}_k(\mathbf{x}, \mathbf{x}_W)}{\partial x_i}  = \\
&\left( \frac{\partial \mathbf{K}^{-1}(\mathbf{x})}{\partial x_i} + \mathbf{K}^{-1}(\mathbf{x})\frac{\partial \mathbf{K}(\mathbf{x}_W)} {\partial x_i} \mathbf{K}^{-1}(\mathbf{x}_W) \right) \mathbf{f}^{eq}_{k}(\mathbf{x}_W) \\
&+ \frac{\partial \mathbf{K}^{-1}(\mathbf{x})}{\partial x_i}  \nabla{\mathbf{f}^{eq}_{k}(\mathbf{x}_W)} \left( \mathbf{x}-\mathbf{x}_W \right) \\
&+\mathbf{K}^{-1}(\mathbf{x}) \mathbf{K}(\mathbf{x}_W) \frac{\partial \mathbf{u}_k(\mathbf{x}_W)} {\partial x_i}.
\end{split}
\end{equation}
Eq.~\eqref{eq:duESL_tilde_dx} and \eqref{eq:duFESLf_tilde_dx_2} will be used to asses the optimality conditions in the next section.

\section{Assessment of the optimality conditions}
\label{subsec:kktcondcomp}
To guarantee the possibility that the proposed F-ESL approach can recognize a feasible local, or even global, optimum of the original dynamic response optimization problem, one must verify that with F-ESL the Karush–Kuhn–Tucker (KKT) optimality conditions of the original dynamic response optimization problem are satisfied at the optimal solution obtained by F-ESL. 

The KKT optimality conditions consist of the stationarity conditions:
\begin{equation}
    \label{eq:KKT_stat_cond_def}
    \nabla \mathcal{L}(\mathbf{x}, \bm{\lambda}) = 
    \nabla f_c(\mathbf{x}) + 
    \sum_{j=1}^{m} \lambda_j \nabla g_j(\mathbf{x}),
\end{equation}
the primal and dual feasibility conditions:
\begin{equation}
    \label{eq:KKT_primal_cond_def}
     \begin{aligned}
        g_j(\mathbf{x}) & \leq 0, \quad j = 1, 2, \ldots, m \\
        \lambda_j & \geq 0, \quad j = 1, 2, \ldots, m,
    \end{aligned}
\end{equation}
and the complementary conditions:
\begin{equation}
    \label{eq:KKT_comp_cond_def}
    \lambda_j g_j(\mathbf{x}) = 0, \quad j = 1, 2, \ldots, m
\end{equation}
where $\mathbf{\lambda}$ is the vector of Lagrange multipliers associated with each constraint considered, including the lower and upper bounds of the optimization variables.

At an optimal point $\mathbf{x}^*$, the KKT stationarity conditions require that the derivatives of the Lagrangian with respect to every design variable converge numerically to zero. 
These values are monitored in the numerical optimization analyses of the examples of Sec.~\ref{sec:numres} to check the proximity of the final solution to an optimal solution\footnote{Even if in a local sense.}. 
When problem \eqref{eq:statoptform} is solved within each sub-problem of the ESL or F-ESL approach, the optimizer yields Lagrange multipliers associated with the current optimal solution of that sub-problem. 
Thus, upon convergence of the ESL or F-ESL approach, the Lagrange multipliers obtained from the final sub-problem solution are then used to compute the stationarity conditions for the original and approximate problem formulations (ESL or F-ESL), respectively, as follows:
\begin{equation}
    \label{eq:stat_cond_ref}
    \nabla \mathcal{L}(\mathbf{x}^{*}, \bm{\lambda}) = 
    \nabla f_c(\mathbf{u}_k(\mathbf{x}^{*}), \mathbf{x}^{*}) + 
    \sum_{j=1}^{m} \lambda_j \nabla g_j(\mathbf{u}_k(\mathbf{x}^{*}),\mathbf{x}^{*})
\end{equation}
and
\begin{equation}
    \label{eq:stat_cond_ESL_FESL}
    \nabla \mathcal{L}(\mathbf{x}^{*}, \tilde{\bm{\lambda}}) = 
    \nabla f_c(\tilde{\mathbf{u}}_k(\mathbf{x}^{*}), \mathbf{x}^{*}) + 
    \sum_{j=1}^{m} \tilde{\lambda}_j \nabla g_j(\tilde{\mathbf{u}}_k(\mathbf{x}^{*}), \mathbf{x}^{*})
\end{equation}
for $k=1,\dots,n_t$. In Eq.~\eqref{eq:stat_cond_ESL_FESL}, the tilde (i.e., $\tilde{ \cdot }$) denotes the values corresponding to either the ESL or F-ESL approach.

In this regard, in \cite{stolpe2014equivalent} Stolpe analytically shows that for dynamic response problems, the stationarity conditions are equivalent between the original and approximate problems if the displacement sensitivities (for the original dynamic problem, $\frac{\partial \mathbf{u}_k(\mathbf{x})}{\partial x_i}$, and its static approximation, $\frac{\partial \tilde{\mathbf{u}}_k(\mathbf{x})}{\partial x_i}$) are identical at the final point $\mathbf{x}^*$. 
This translates into requiring that the approximate approach adopted is effectively a first-order approximation of the original problem. 
In the analytical derivation presented here, we strengthen the stopping criteria by setting $\epsilon=0$. Let us compute the displacement sensitivities at convergence, where $\tilde{\mathbf{x}}=\mathbf{x}^*$, for the three discussed formulations.
For ESL, using Eq.~\eqref{eq:duESL_tilde_dx} at the optimal point:
\begin{equation}
    \label{eq:duESL_tilde_dx_atConv}
    \begin{split}
    \frac{\partial \tilde{\mathbf{u}}^{ESL}_k(\mathbf{x}^*)}{\partial x_i} = \frac{\partial \mathbf{K}^{-1}(\mathbf{x}^*)}{\partial x_i} \mathbf{f}^{ESL}_{k}(\mathbf{x}^*).
    \end{split}
\end{equation}
Using Eq.~\eqref{eq:duESL_tilde_dx_atConv} and the definition in Eq.~\eqref{eq:du_dx}, we obtain:
\begin{equation}
    \label{eq:duESL_tilde_dx_atConv_2}
    \begin{split}
    \frac{\partial \tilde{\mathbf{u}}^{ESL}_k(\mathbf{x}^*)}{\partial x_i} = 
    \frac{\partial \mathbf{u}_k(\mathbf{x}^*)}{dx_i} -\mathbf{K}^{-1}(\mathbf{x}^*) \frac{\partial \mathbf{f}^{eq}_{k}(\mathbf{x}^*)}{\partial x_i}.
    \end{split}
\end{equation}
For the F-ESL approach proposed here, using Eq.~\eqref{eq:duFESLf_tilde_dx_2} at convergence and the following result
\begin{equation}
\mathbf{K}^{-1}(\mathbf{x}^*)\frac{\partial \mathbf{K}(\mathbf{x}^*)} {\partial x_i} \mathbf{K}^{-1}(\mathbf{x}^*)=-\frac{\partial \mathbf{K}^{-1}(\mathbf{x}^*)} {\partial x_i},
\end{equation}
we obtain the following:
\begin{equation}
\label{eq:duFESLf_tilde_dx_atConv}
\begin{split}
&\frac{\partial \tilde{\mathbf{u}}^{FESL}_k(\mathbf{x}^*)}{\partial x_i}  = 
\left( \frac{\partial \mathbf{K}^{-1}(\mathbf{x}^*)}{\partial x_i} -\frac{\partial \mathbf{K}^{-1}(\mathbf{x}^*)} {\partial x_i} \right) \mathbf{f}^{eq}_{k}(\mathbf{x}^*) +\\
& \frac{\partial \mathbf{K}^{-1}(\mathbf{x}^*)}{\partial x_i}  \nabla{\mathbf{f}^{eq}_{k}(\mathbf{x}^*)} \left( \mathbf{x}^*-\mathbf{x}^* \right)+ 
\mathbf{K}^{-1}(\mathbf{x}^*) \mathbf{K}(\mathbf{x}^*) \frac{\partial \mathbf{u}_k(\mathbf{x}^*)} {\partial x_i}.
\end{split}
\end{equation}
At the end of the optimization process, where the stopping condition in Eq.~\eqref{eq:stoppingcond} is satisfied, the first two terms equal zero, and therefore Eq.~\eqref{eq:duFESLf_tilde_dx_atConv} reduces to:
\begin{equation}
\label{eq:duFESLf_tilde_dx_atConv2}
\begin{split}
\frac{\partial \tilde{\mathbf{u}}^{FESL}_k(\mathbf{x}^*)}{\partial x_i}  =  \frac{\partial \mathbf{u}_k(\mathbf{x}^*)}{\partial x_i}.
\end{split}
\end{equation}

From these analytical derivations it becomes clear that the displacement sensitivities of ESL shown in Eq.~\eqref{eq:duESL_tilde_dx_atConv_2} differ from the displacement sensitivities of the original dynamic response optimization problem at the final optimal point. 
Thus, ESL can not detect an optimal solution for the original problem. 
On the other hand, the displacement sensitivities of the F-ESL formulation, shown in Eq.~\eqref{eq:duFESLf_tilde_dx_atConv2}, approach those of the original problem as the F-ESL algorithm converges towards the final optimal solution.

\section{Numerical examples}
\label{sec:numres}

The basic ESL and the proposed first-order F-ESL formulations were implemented in MATLAB R2023a and applied to a suite of structural optimization problems with different constraints and dynamic loading. 
In the examples, the structural dynamic response analyses were performed using the average acceleration scheme of the Newmark-$\beta$ algorithm. The optimization problems were solved with the Sequential Quadratic Programming (SQP) algorithm implemented in the \verb+fmincon+ MATLAB function with its default settings, unless otherwise specified. In order to compute the displacement sensitivities in each of the examples, we relied on the Direct Sensitivity Analysis method, and more details on this can be found in \ref{subsec:directsen_andx}.

\subsection{An illustrative example}
\label{subsec:studycasedes}

\begin{figure}[htbp!]
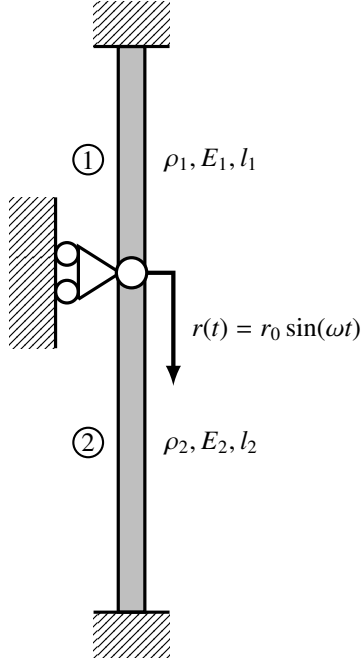

    \centering
    \include{elements/system_drawing_example1}
    \setlength{\abovecaptionskip}{-0.3cm}
    \setlength{\belowcaptionskip}{-0.2cm}
    \caption{\small Example of Sec.~\ref{subsec:studycasedes}. A two-bar structure with one degree of freedom subjected to a vertical sinusoidal force.}
    \label{fig:systemDrawing_example1}
\end{figure}

The first example considered in this paper consists of a simple structure with one degree of freedom, similar to the one discussed by Stolpe et al. \citep{stolpe2018equivalent}. The structure consists of two aligned bars with identical material properties but with different lengths, as shown in Fig.~\ref{fig:systemDrawing_example1}. A longitudinal sinusoidal force acts between the two bars, causing axial displacements in both bars. There are two design variables, $x_1$ and $x_2$, which are the cross-sectional areas of the bars. The structure has one degree of freedom, consisting of the vertical displacement (and associated velocity and acceleration) of the internal hinge to which the load is applied. 

The optimization problem considered here aims to minimize the structural compliance, i.e. $\int_{t_0}^{t_f}r(t) u(t) dt = \int_{t_0}^{t_f}u(t) k\, u(t)dt$, with a volume constraint. The problem was intentionally chosen in \cite{stolpe2018equivalent} to study specific characteristics of the original zeroth-order ESL approach. Note that the parameters' values in this example were chosen arbitrarily to show the behavior of the problem, and therefore will be presented with no specified units of measurement. 
The dynamic response optimization problem and the static response optimization approximated sub-problems are convex problems, and therefore can be solved, allowing for a comprehensive investigation. The static sub-problems are ensured to have non-empty feasible sets. Additionally, because there are only two design variables, it is possible to provide a graphic visualization of the design domain and numerical behavior of the optimization approaches discussed. 

The reference dynamic response optimization problem formulation reads:
\begin{equation}
\label{eq:dynoptform_example1}
\begin{split}
    \underset{\mathbf{x} \in {\rm I\!R}^2}{\text{minimize }} &  f_c(\mathbf{x},u_k (\mathbf{x}))=\Delta t \sum_{k=s_t}^{n_t}{k (\mathbf{x})} u_k^2 (\mathbf{x})\\
    \text{subject to:  } & g(\mathbf{x})= \sum_{i=1}^{2}{\frac{x_il_i}{V_{max}}}-1 \leq 0 \\
    & x_i^{min} \leq x_i \leq x_i^{max}, \hspace{5mm} i=1,2 \\
    \text{with:  } &m(\mathbf{x}) \ddot{u}_k(\mathbf{x})
    +c \dot{u}_k(\mathbf{x}) 
    +k(\mathbf{x}) u_k(\mathbf{x}) = r_0 \sin(\omega t_k)\\
    & \text{for } k=1,2,\dots,n_t
\end{split}
\tag{$\mathcal{P}^1_D$}
\end{equation} 
where $V_{max}=1$ is the maximum allowed volume for the entire structure, and the design variables, i.e. the bars' cross-sectional areas, are allowed to vary between $x_i^{min}=0.1$ and $x_i^{max}=1$. 
With reference to Fig.~\ref{fig:systemDrawing_example1}, the modulus of elasticity of each bar is $E_1=E_2=1$, the densities are $\rho_1=\rho_2=1$, and the lengths of the bars are $l_1=1$ and $l_2=3/2$. The representative lumped-mass and stiffness expressions related to the single DOF considered become
\begin{equation}
    \label{eq:stiff_example1}
    k(\mathbf{x}) = \frac{x_1 E_1}{l_1} + \frac{x_2 E_2}{l_2} = x_1 + \frac{2}{3}x_2,
\end{equation}
and
\begin{equation}
    \label{eq:mass_example1}
    m(\mathbf{x}) = \frac{1}{2} x_1 \rho_1 l_1 + \frac{1}{2} x_2 \rho_2 l_2 = \frac{1}{2} x_1 + \frac{3}{4} x_2.
\end{equation}
The damping for the structure at the free degree of freedom was chosen to be independent of the mass and stiffness values, with a constant value of $c=0.1$. The sinusoidal dynamic load amplitude and frequency are $r_0=1$ and $\omega=\pi/2$, respectively.

The ESL and F-ESL optimization problem formulations at the $W$-th iteration, following the definitions in \eqref{eq:statoptform}, are: 
\begin{equation}
\label{eq:statoptform_example1}
\begin{split}
    \underset{\mathbf{x} \in {\rm I\!R}^2}{\text{minimize }} &  f_c(\mathbf{x},\tilde{u}_k (\mathbf{x})) = \Delta t \sum_{k=s_t}^{n_t}{k (\mathbf{x})} \tilde{u}_k^2 (\mathbf{x})\\
    \text{subject to:  } & g(\mathbf{x})= \sum_{i=1}^{n}{\frac{x_il_i}{V_{max}}}-1 \leq 0 \\
    & x_i^{min} \leq x_i \leq x_i^{max}, \hspace{5mm} i=1,2 \\
    \text{with:  } & k (\mathbf{x}) \tilde{u}_k(\mathbf{x}, \mathbf{x}_W) 
     = f^{ESL}_{k} (\mathbf{x}_W)\\ 
    \text{or } & k (\mathbf{x}) \tilde{u}_k(\mathbf{x}, \mathbf{x}_W) 
     = f^{FESL}_{k} (\mathbf{x}, \mathbf{x}_W)\\ 
     &\text{for } k=s_t,\dots,n_t.
\end{split}
\tag{$\mathcal{P}^1_{S,W}$}
\end{equation}
In the numerical results reported in this section, as in \cite{stolpe2018equivalent}, the time step value taken for numerical integration is $\Delta t=0.2$, and the dynamic analysis is performed in the time span of $t\in[0,t_1+\pi/\omega]$. 
The objective function is computed over the time window $t\in[t_1, t_1+\pi/\omega]$ with $t_1=200$ and $\omega=\pi/2$, that is, $t\in[200, 202]$.
The optimization variables are both initialized with the same value $x_1^0=x_2^0=0.2$, which is within the feasible domain $[0.1,1]$.

The structural system considered (Fig.~\ref{fig:systemDrawing_example1}) was first optimized solving directly the original dynamic response optimization problem formulation \eqref{eq:dynoptform_example1}, and then with the ESL and F-ESL approaches solving a sequence of sub-problems \eqref{eq:statoptform_example1} until convergence, with the stopping criteria value of $\epsilon=1\mathrm{e}{-10}$ in Eq.~\eqref{eq:stoppingcond}. 
Table~\ref{tab:resultstable_example1} summarizes the results obtained with ESL, F-ESL, and by solving the original dynamic optimization problem formulation directly.
The table lists the number of time-history analyses performed; 
the final design obtained; 
the final value of the objective function; 
the gradient of the objective function in correspondence of the optimal solution; 
the constraint value; 
the gradient of the Lagrangian function of the problem in correspondence of the optimal solution (see Sec.~\ref{subsec:kktcondcomp}).

\begin{figure}[htbp!]
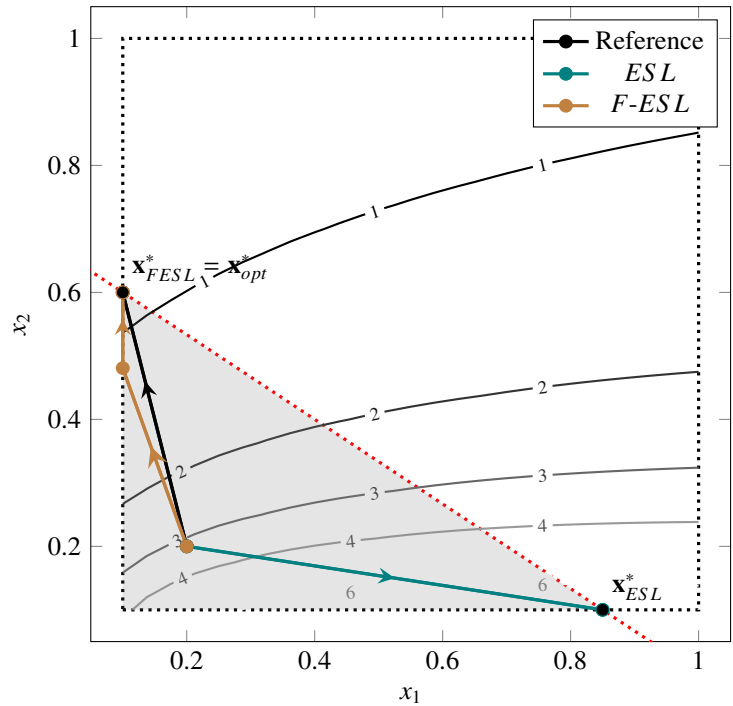

    \centering
    \include{elements/VariablePath_plot_example1}
    \setlength{\abovecaptionskip}{-0.5cm}
    \setlength{\belowcaptionskip}{-0.1cm}
    \caption{\small Illustration of the numerical example of Sec.~\ref{subsec:studycasedes}. The figure shows the design domain, the contour of the objective function of the original dynamic problem, and the optimization paths across the iterations for ESL, F-ESL and the solution of the reference (original) problem. The area in gray defines the feasible design domain, enclosed by the box constraints and the volume constraints, marked by black and red dotted lines, respectively.}
    \label{fig:FESLplot}
\end{figure}

The ESL algorithm terminated after $2$ outer iterations, and the F-ESL algorithm after $3$ outer iterations, where at each outer iteration a linear dynamic response analysis is performed. For comparison, the direct solution of the reference dynamic problem required in total $7$ time-history analyses to converge.
It should be mentioned, that in this example, during the optimization analyses with ESL and F-ESL moving limits for the design variables were not considered.
As shown in Table.~\ref{tab:resultstable_example1}, ESL converged to a design that is not an optimal solution of the original problem. F-ESL, instead, converges to the same optimal solution as that obtained by solving the original problem. 
\begin{table*} [htbp!]\setlength{\tabcolsep}{4pt}
    \caption{\small Numerical example of Sec. \ref{subsec:studycasedes}. Optimization results obtained with ESL, F-ESL, and by solving the original problem directly. The table provides the optimal design $\mathbf{x}^*$, the number of time-history analyses required $N_{THA}$, the value of the objective function $f_c(\mathbf{x}^*)$, the objective function gradient $\nabla f_c(\mathbf{x}^*)$, the constraint value $g(\mathbf{x}^*)$, and the gradient of the Lagrangian of the problem $\nabla \mathcal{L}(\mathbf{x}^*)$, all evaluated at the final optimal solution $\mathbf{x}^*$.}
    \centering \scriptsize
    \include{elements/results_table_example1}
    \label{tab:resultstable_example1}
\end{table*}
As reported also in \cite{stolpe2018equivalent}, the optimal design for \eqref{eq:dynoptform_example1} is $\mathbf{x}^*_{opt} =
    \begin{bmatrix}
        0.1 &
        0.6
    \end{bmatrix}^{T}$.
Because this numerical example consists of only two design variables, the trajectory of the design variables for each approach during the optimization process can be visualized in a 2D plot, as shown in Fig.~\ref{fig:FESLplot}.
An additional optimization analysis was performed with ESL and F-ESL, this time using the known optimal solution $\mathbf{x}^*_{opt}$ directly as the starting point. Even in this case, the ESL algorithm failed to recognize the optimal solution and converged to the same worst final design solution $\mathbf{x}^*_{ESL} \neq \mathbf{x}^*_{opt}$. This highlights the inability of ESL to recognize optimal solutions of the original problem. This limitation does not seem to depend on the starting point or its proximity to an optimal solution. 

\begin{figure*}
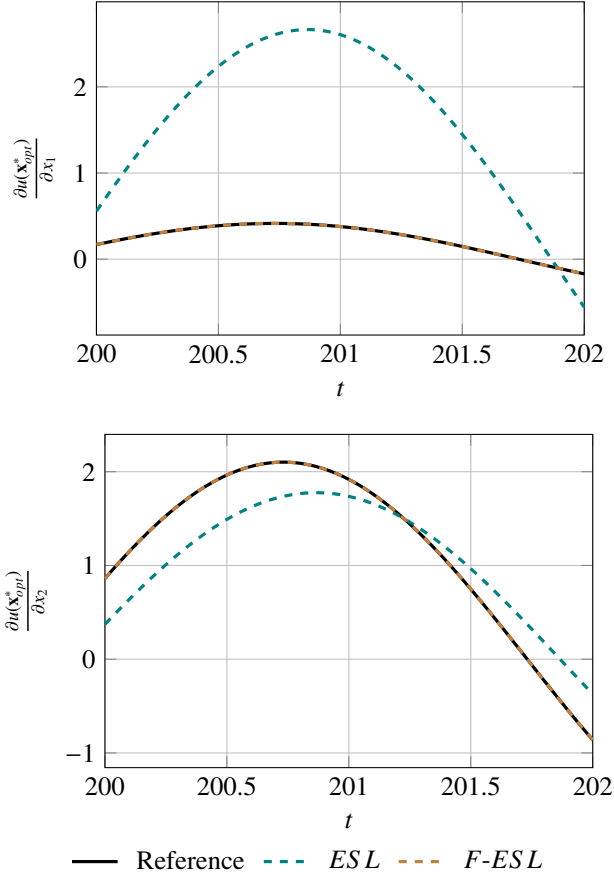

    \centering
    \include{elements/dUdx_xfesl_plot}
    \setlength{\abovecaptionskip}{-0.5cm}
    \setlength{\belowcaptionskip}{-0.1cm}
    \caption{\small Numerical example of Sec. \ref{subsec:studycasedes}. Comparison of the displacements sensitivities with respect to the design variables $x_1$ and $x_2$, of the original dynamic-response optimization problem, the ESL and the F-ESL formulations, at the final solution of the F-ESL, which is the optimal solution of the dynamic-response optimization problem \eqref{eq:dynoptform_example1}.}
    \label{fig:dUdxplot_xfesl}
\end{figure*}

\begin{figure}
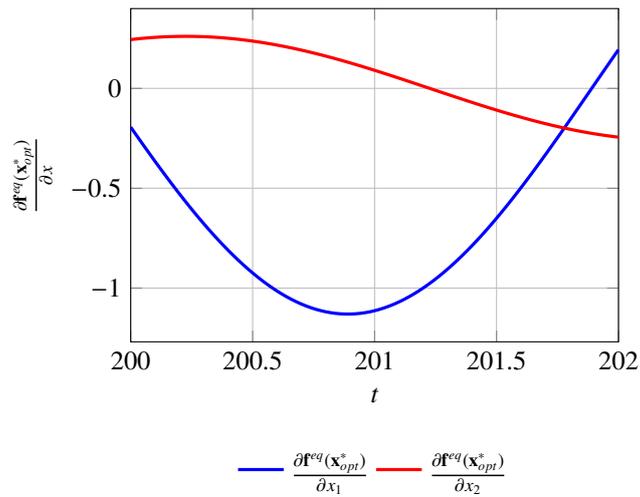

    \centering
    \include{elements/dFesl_dA}
    \setlength{\abovecaptionskip}{-0.5cm}
    \setlength{\belowcaptionskip}{-0.1cm}
    \caption{\small Numerical example of Sec. \ref{subsec:studycasedes}. Gradient of the equivalent static loads (Eq.~\eqref{eq:nablafESL}) at the optimal solution $\mathbf{x}^*_{opt}$ of the dynamic problem \eqref{eq:dynoptform_example1}.
    }
    \label{fig:dFesldAplot}
\end{figure}

The main crucial difference between the ESL and the F-ESL algorithms is that the F-ESL formulation is defined mathematically so that when it converges to the final optimal solution, the displacement sensitivities of the static problem are identical to those of the original dynamic response problem, as discussed in Sec.~\ref{subsec:kktcondcomp}. 
In support of this statement, the displacement sensitivities of ESL, F-ESL, and the original dynamic response problem are shown in Fig.~\ref{fig:dUdxplot_xfesl}, for the two design variables $x_1$ and $x_2$. The sensitivities are computed at the optimal solution of the original dynamic response optimization problem, $\mathbf{x}^*_{opt}$, to which the F-ESL approach successfully converges.
\begin{figure*}
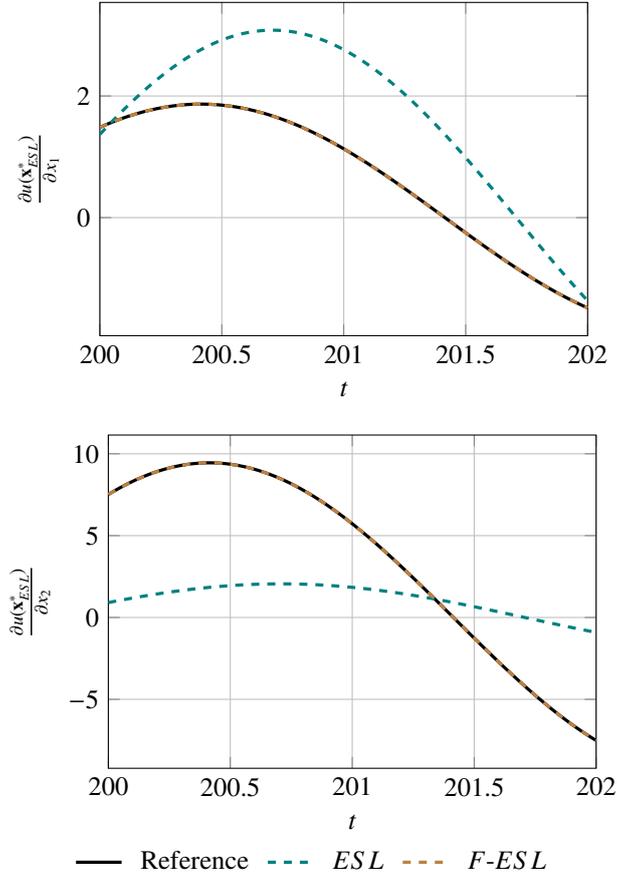

    \centering
    \include{elements/dUdx_xesl_plot}
    \setlength{\abovecaptionskip}{-0.5cm}
    \setlength{\belowcaptionskip}{-0.1cm}
    \caption{\small Numerical example of Sec. \ref{subsec:studycasedes}. Comparison of the displacements sensitivities with respect to the design variables $x_1$ and $x_2$, of the original dynamic-response optimization problem, the ESL and the F-ESL formulations, at the final solution of the ESL, which is \textit{not} the optimal solution of the dynamic-response optimization problem \eqref{eq:dynoptform_example1}.}
    \label{fig:dUdxplot_xesl}
\end{figure*}
The displacement sensitivities with respect to each variable of F-ESL coincide with those of the original dynamic response problem at $\mathbf{x}=\mathbf{x}^*_{opt}$. The displacement sensitivities values of ESL, instead, are different from those of the original problem, as shown also analytically in Eq.~\eqref{eq:duESL_tilde_dx_atConv_2}. 
The same behavior is observed in Fig.~\ref{fig:dUdxplot_xesl}, at the final solution point of the ESL algorithm. 
In addition, the gradient components of the equivalent static loads $\mathbf{f}_k^{eq}(\mathbf{x})$ at the optimal design $\mathbf{x}^*_{opt}$ with respect to $x_1$ and $x_2$ are plotted in Fig.~\ref{fig:dFesldAplot}.  
It can be clearly observed that, as we discussed previously in Sec.~\ref{subsec:kktcondcomp}, the gradient of the equivalent static loads is not, in general, zero in correspondence with an optimal solution of the problem.

\subsection{Seismic design optimization of a two-story shear frame}
\label{subsec:secondex}
\begin{figure}[htbp!]
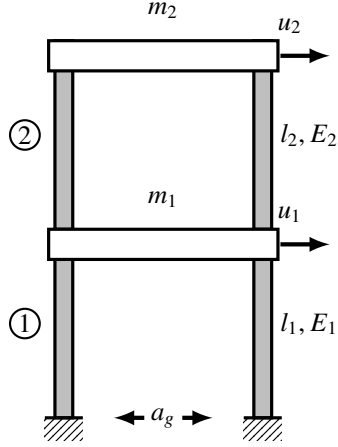

    \centering
    \include{elements/system_drawing_example2}
    \setlength{\abovecaptionskip}{-0.3cm}
    \setlength{\belowcaptionskip}{-0.2cm}
    \caption{\small Numerical example of Sec. \ref{subsec:secondex}. Two-story shear frame structure subjected to ground excitation $a_g(t)$ with two degrees of freedom, $u_1$ and $u_2$.}
    \label{fig:systemDrawing_example2}
\end{figure}

In this section, we present and discuss an additional numerical example. It concerns the design optimization of a shear frame subjected to seismic excitation (i.e., ground acceleration due to an earthquake). This example is inspired by a similar example discussed in \cite{idels2020performance}.
In particular, the structure considered consists of a two-story shear frame with two degrees of freedom $u_1$ and $u_2$. 
The structure is subjected to seismic excitation at the ground level, as shown in Fig.~\ref{fig:systemDrawing_example2}. The masses of the first and second floors are $m_1=m_2=4$ ton. Each floor is made of two identical steel columns. Assuming that the columns have a circular profile, the optimization variables are the columns' cross-section diameters $x_1$ and $x_2$. Moreover, the length of each column is $l_1=l_2=3.5$ m and the Young's modulus considered is $E_1=E_2=200$ GPa. Both floors are considered rigid horizontally, and as a consequence, both columns of each floor experience the same horizontal motion during the earthquake excitation.

Given the earthquake ground acceleration $a_{g,k}$ at each time step $k$, the resulting force applied to each mass is:
\begin{equation}
    \label{eq:earthquakeforce}
    \mathbf{f}_{g,k} = - \mathbf{M} \mathbf{e} a_{g,k}.
\end{equation}
In Eq.~\eqref{eq:earthquakeforce} $\mathbf{M}$ is the global mass matrix of the system and $\mathbf{e}$ is the influence vector that applies the ground acceleration to each mass (its entries are all equal to one).
The stiffness of each column associated to an horizontal displacement of one of its ends is:
\begin{equation}
    \label{eq:columnStiffness}
    k_i=k(x_i)=\frac{3 E_i I_i(x_i)}{L_i^3}
\end{equation}
for $i=1,2$. In addition, in Eq.~\eqref{eq:columnStiffness} $I(x_i)$ is the columns' moment of inertia which for a circular cross section is defined as follows:
\begin{equation}
    \label{eq:inertia_term}
    I(x_i)=\frac{\pi x_i^4}{64}.
\end{equation} 
The stiffness and mass matrices and the influence vector are therefore defined as follows: 
\begin{equation}
\begin{split}
    \label{eq:stiff_mass_matrices}
    &\mathbf{K}(\mathbf{x})=\begin{bmatrix}
        2k_1+2k_2  & -2k_2 \\
        -2k_2 & 2k_2 \\
    \end{bmatrix}, \; 
    \mathbf{M}=\begin{bmatrix}
        m_1 & 0 \\
        0 & m_2 \\
    \end{bmatrix}, \\ 
    & \mathbf{e}=\begin{bmatrix}
        1  \\
        1 \\
    \end{bmatrix}.
\end{split}
\end{equation}
In the previous equation $k_1$ and $k_2$ are multiplied by two because each floor is made of two identical columns. 
We include the inherent structural damping by means of the Rayleigh damping matrix formulation \citep{hudson1995dynamics}:
\begin{equation}
    \label{eq:damping_matrix}
    \mathbf{C}_s(\mathbf{x}) = \alpha_M(\mathbf{x}) \mathbf{M} + \alpha_K(\mathbf{x}) \mathbf{K}(\mathbf{x}).
\end{equation}
The Rayleigh coefficients, $\alpha_M(\mathbf{x})$ and $\alpha_K(\mathbf{x})$, depend on the two natural frequencies of the structure $\omega_1$ and $\omega_2$, computed by solving the following eigenvalue problem:
\begin{equation}
    \label{eq:eigenvalueprob}
    (\mathbf{K}-\omega^2_i\mathbf{M})\bm{\phi}_i=\mathbf{0}
\end{equation}
in which $\bm{\phi}_i$ is the eigenvector corresponding to each natural frequency $\omega_i$. In practice, the eigenvalues $\omega_i^2$ are computed by passing to the MATLAB \verb+eig+ function the matrix $\mathbf{M}^{-1}\mathbf{K}$.
For each mode $i$, the damping ratio is defined by the expression:
\begin{equation}
    \label{eq:dampingratioeq}
    \xi_i=\frac{1}{2}(\frac{\alpha_M}{\omega_i}+\alpha_K\omega_i).
\end{equation}
If we consider a damping ratio $\xi=\xi_1=\xi_2=5\%$, the Rayleigh damping coefficients can be calculated by solving the linear system of equations of Eq.~\eqref{eq:dampingratioeq}, leading to the following result:
\begin{equation}
    \label{eq:rayleighcoeff}
    \alpha_M=2\xi\frac{\omega_1 \omega_2}{\omega_1 + \omega_2}, \quad \alpha_K=2\xi\frac{1}{\omega_1 + \omega_2}.
\end{equation}

\begin{figure}[htbp!]
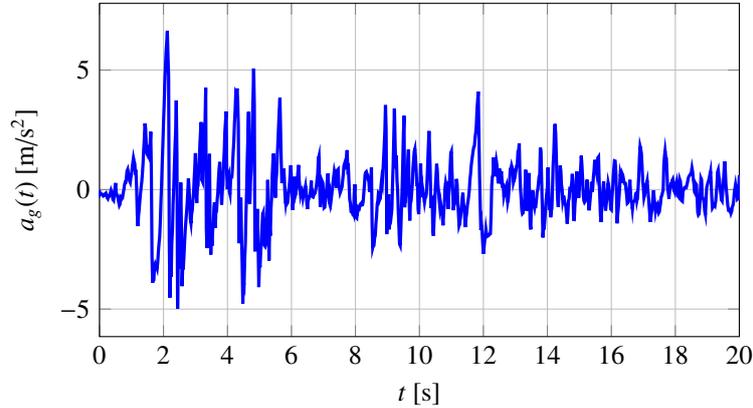

    \centering
    \include{elements/earthquake_plot}
    \setlength{\abovecaptionskip}{-0.3cm}
    \setlength{\belowcaptionskip}{-0.2cm}
    \caption{\small Numerical example of Sec. \ref{subsec:secondex}. Plot of the first $20$ sec of the earthquake ground acceleration record ($a_g(t)$) from the “LA $10$\% in $50$ years” ground motions ensemble \citep{somerville1997development}.}
    \label{fig:earthquakeplot}
\end{figure}

The structural volume of the shear frame considered in this example is minimized, while considering a constraint of the maximum inter-story drift\footnote{The inter-story drift is the relative horizontal displacement of the columns' ends.} value in time. The drift vector $\mathbf{d}_k$ in each $k$-th time step is computed by means of a linear transformation of the displacement vector $\mathbf{u}_k$:
\begin{equation}
    \label{eq:driftcalc}
    \mathbf{d}_k=\mathbf{D}\mathbf{u}_k, \quad \mathbf{D}=\begin{bmatrix}
        1 & 0 \\
        -1 & 1 \\
    \end{bmatrix}.
\end{equation} 
Thus, the optimization problem considered in this example is the following:
\begin{equation}
\label{eq:dynoptform_example2}
\begin{split}
    \underset{\mathbf{x} \in {\rm I\!R}^2}{\text{minimize }} &  f_c(\mathbf{x})=\sum_{i=1}^{2}{2l_i \frac{\pi x_i^2}{4}}\\
    \text{subject to:  } & g_{i,k}(\mathbf{x},\mathbf{u}_k(\mathbf{x}))= \frac{|d_i(\mathbf{u}_k(\mathbf{x}))|}{d_{max}} - 1 \leq 0 \\
    & x_i^{min} \leq x_i \leq x_i^{max}, \hspace{5mm} i=1,2 \quad k=1,2,\dots,n_t \\
    \text{with:  } &\mathbf{M} \mathbf{\ddot{u}}_k(\mathbf{x})
    +\mathbf{C}_s(\mathbf{x}) \mathbf{\dot{u}}_k(\mathbf{x}) 
    +\mathbf{K}(\mathbf{x}) \mathbf{u}_k(\mathbf{x}) = \mathbf{f}_{g,k}\\
    & \text{for } k=1,2,\dots,n_t.
\end{split}
\tag{$\mathcal{P}^{2}_D$}
\end{equation}
In \eqref{eq:dynoptform_example2} the design variables can vary between $x_{min}=0.001$ m and $x_{max}=0.5$ m, and the maximum allowable inter-story drift is $d_{max}=0.1$ m. As mentioned previously, we considered the first $n_t=20$ s of the LA02 ground acceleration record, from the “LA $10$\% in $50$ years” ground motions ensemble \citep{somerville1997development}, with a constant time step $\Delta t=0.02$ s. The ground acceleration record considered is shown in Fig.~\ref{fig:earthquakeplot}. As a result in the optimization problem considered here, there are in total $2002$ constraints. 

In this example, we focus only on the use of F-ESL for solving the dynamic response optimization problem at hand. 
In each $W$-th iteration of the F-ESL approach, the static response optimization sub-problem is formulated as follows:
\begin{equation}
\label{eq:statoptform_example2}
\begin{split}
    \underset{\mathbf{x} \in {\rm I\!R}^2}{\text{minimize }} &  f_c(\mathbf{x})=\sum_{i=1}^{2}{2l_i \frac{\pi x_i^2}{4}}\\
    \text{subject to:  } & g_{i,k}(\mathbf{x},\mathbf{\tilde{u}}_k(\mathbf{x}))= \frac{|d_i(\mathbf{\tilde{u}}_k(\mathbf{x}))|}{d_{max}} - 1 \leq 0 \\
    & x_i^{min} \leq x_i \leq x_i^{max}, \hspace{5mm} i=1,2 \quad k=1,2,\dots,n_t \\
    \text{with:  } &\mathbf{K} (\mathbf{x}) \mathbf{\tilde{u}}_k(\mathbf{x}, \mathbf{x}_W) 
     = \mathbf{f}^{FESL}_{k} (\mathbf{x}, \mathbf{x}_W)\\
    & \text{for } k=1,2,\dots,n_t.
\end{split}
\tag{$\mathcal{P}^2_{S,W}$}
\end{equation}
The equivalent static loads considered in the F-ESL formulation and their derivatives with respect to the design variable are calculated using Eq.~\eqref{eq:feslfdef} and Eq.~\eqref{eq:nablafESL} (see Algorithm~\ref{pseudoalgo}) with $\epsilon=1\mathrm{e}{-6}$ for the stopping criteria defined in Eq.~\eqref{eq:stoppingcond}.
In the \ref{subsec:driftconstsen_andx}, the details of the sensitivity analysis of the drift constraint are provided.
For comparison purposes, the original dynamic response problem \eqref{eq:dynoptform_example2} was also solved directly using the gradient-based SQP algorithm of the MATLAB \verb+fmincon+ function and a gradient-free genetic algorithm optimizer (GA), provided by the MATLAB \verb+ga+ function. The \verb+fmincon+ and \verb+ga+ functions were used with default settings, except that with 
the SQP solver of \verb+fmincon+  \verb+ConstraintTolerance+ was set to $=1e-5$ and \verb+StepTolerance+ was set to $1e-12$; with 
\verb+ga+ the constraint tolerance violation was set to \verb+ConstraintTolerance=1e-5+. 
The initial value for all design variables, that is, the initial diameters of the beams, is $x_i^0=0.5$ m. 
 
\begin{figure}[htbp!]
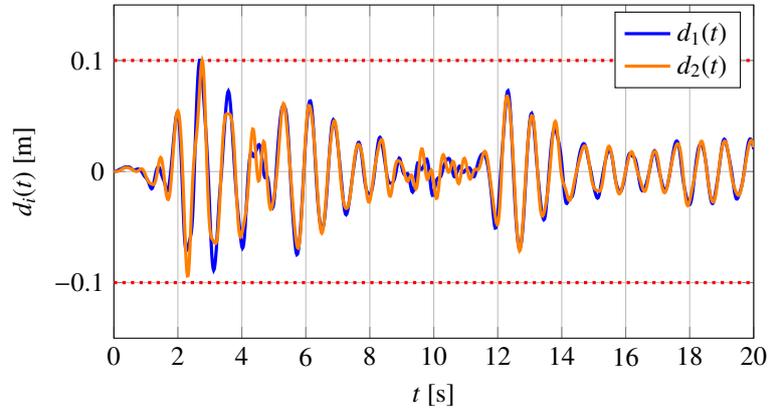

    \centering
    \include{elements/drifts_plot_example2}
    \setlength{\abovecaptionskip}{-0.3cm}
    \setlength{\belowcaptionskip}{-0.2cm}
    \caption{\small Numerical example of Sec. \ref{subsec:secondex}. Plot of the inter-story drifts, $d_1(t)=u_1(t)$ and $d_2(t)=u_2(t)-u_1(t)$, of the optimized structure. The maximum and minimum allowed drift values are marked with red dotted lines.}
    \label{fig:driftsplot}
\end{figure}

\begin{figure}[htbp!]
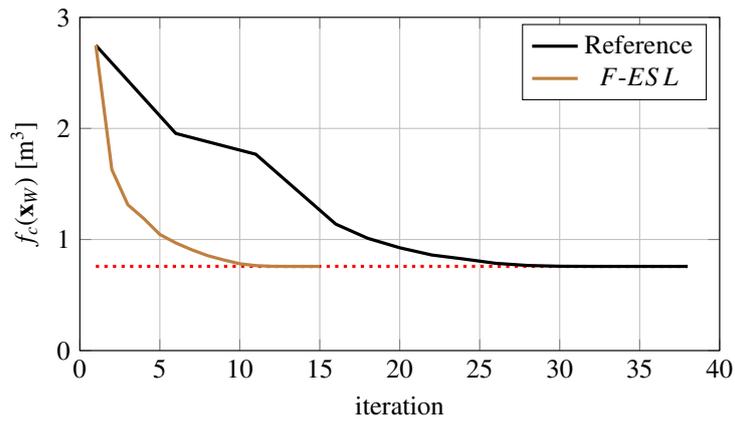

    \centering
    \include{elements/obj_iterations_plot_example2}
    \setlength{\abovecaptionskip}{-0.3cm}
    \setlength{\belowcaptionskip}{-0.2cm}
    \caption{\small Numerical example of Sec. \ref{subsec:secondex}. Objective function values tracked at each iteration (i.e., after every transient analysis) associated with the solution of the original dynamic response problem with SQP, and the proposed F-ESL approach. }
    \label{fig:objiters_example2}
\end{figure}

\begin{figure}[htbp!]
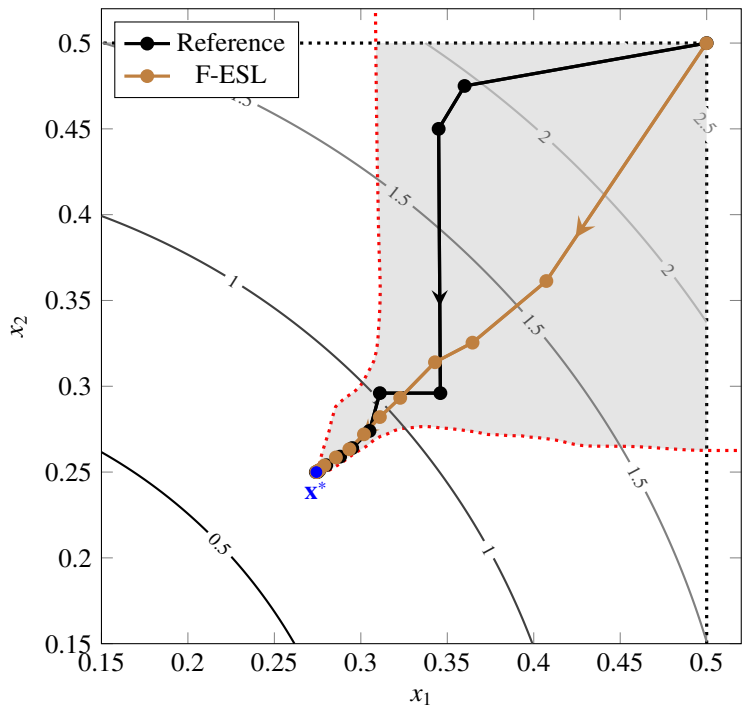

    \centering
    \include{elements/VariablePath_plot_example2}
    \setlength{\abovecaptionskip}{-0.5cm}
    \setlength{\belowcaptionskip}{-0.1cm}
    \caption{\small Numerical example of Sec. \ref{subsec:secondex}. Path of the design points across iterations followed by the F-ESL approach, and by SQP while solving the original problem formulation. The contour shows the compliance function values of the original dynamic response optimization problem. The red dotted line marks the contour $g=0$ that limits the feasible design domain, which is colored in gray.}
    \label{fig:FESLplot_example2}
\end{figure}

Table~\ref{tab:resultstable_example2} gathers, for each optimization approach, the number of time-history analyses performed, the final optimized design, the objective function and its gradients, the maximum constraint value, and the gradient of the Lagrangian. Fig.~\ref{fig:FESLplot_example2} shows the trajectories of the design variables followed by the SQP algorithm while solving the reference dynamic response problem directly, and by the F-ESL algorithm, superimposed on the objective and constraint functions' contours within the design space. 
For the sake of clarity, the plot shows only a selected region of the design domain. 
Looking at Table~\ref{tab:resultstable_example2} as well as Fig.~\ref{fig:FESLplot_example2}, it can be observed that F-ESL converges to the same final design obtained by solving the original dynamic response problem directly with the SQP algorithm (\verb+fmincon)+ and the GA algorithm (\verb+ga+). As expected, GA required a larger number of transient analyses, compared to SQP and F-ESL, to converge to the final solution. 
A direct comparison of computational effort, in terms of the amount of transient analyses performed, highlights the advantage of the F-ESL approach: the static subproblem sequence in \eqref{eq:statoptform_example2} converged after only $15$ transient response evaluations, whereas solving the full dynamic response optimization problem \eqref{eq:dynoptform_example2} with SQP required $38$ transient analyses. 
Table \ref{tab:resultstable_example2} further shows that with F-ESL the gradient components of the Lagrangian at convergence are of the order of $10^{-6}$. This shows that the obtained design meets the optimality conditions of the original problem.  
The inter-story drifts corresponding to the final optimized design obtained with the F-ESL approach are plotted in Fig.~\ref{fig:driftsplot}. It is possible to observe that the optimized response is within the maximum and minimum allowed values of inter-story drifts, marked by red dotted lines. 
In addition, Fig.~\ref{fig:objiters_example2} shows the values of the objective function at each outer iteration (that is, each time a time-history analysis is performed) obtained by solving the original problem directly with SQP and with the proposed F-ESL approach. This comparison further shows that, in this example as well, F-ESL requires less time-history analyses than the direct solution of the original dynamic response optimization problem.  

\begin{table*} [htbp!]\setlength{\tabcolsep}{4pt}
    \caption{\small Numerical example of Sec. \ref{subsec:secondex}. Final results obtained with F-ESL, and by solving the original problem directly with GA and SQP. The table provides the number of time-history-analyses performed $N_{THA}$, the final design $\mathbf{x}^*$, the final value of the objective function $f_c(\mathbf{x}^*)$ its gradient $\nabla f_c(\mathbf{x}^*)$, the gradient of the Lagrangian of the problem $\nabla \mathcal{L}(\mathbf{x}^*)$, the maximum drift constraint $g_{\max}=\max{(g_{i,k}(\mathbf{x^*}))}$.}
    \centering \scriptsize
    \include{elements/results_table_example2}
    \label{tab:resultstable_example2}
\end{table*}

\subsection{Optimization of a truss structure with transient loading and stress constraints}
\label{subsec:thirdex} 
\begin{figure*}
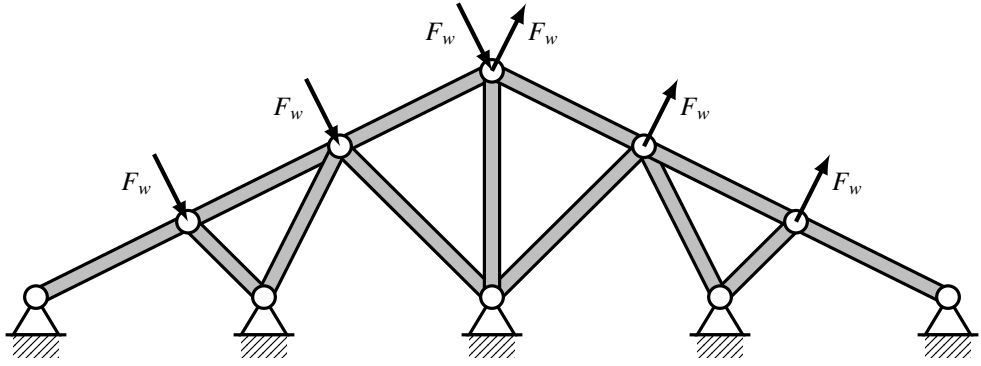
 
    \centering
    \include{elements/system_drawing_example3}
    \setlength{\abovecaptionskip}{-0.3cm}
    \setlength{\belowcaptionskip}{-0.2cm}
    \caption{\small Numerical example of Sec. \ref{subsec:thirdex}. A truss structure with 13 bars, subject to multiple loading conditions.$F_w(t)=p_i(t)$ where $i$ is the load-case considered each time.}
    \label{fig:systemDrawing_example3}
\end{figure*}

\begin{figure*}
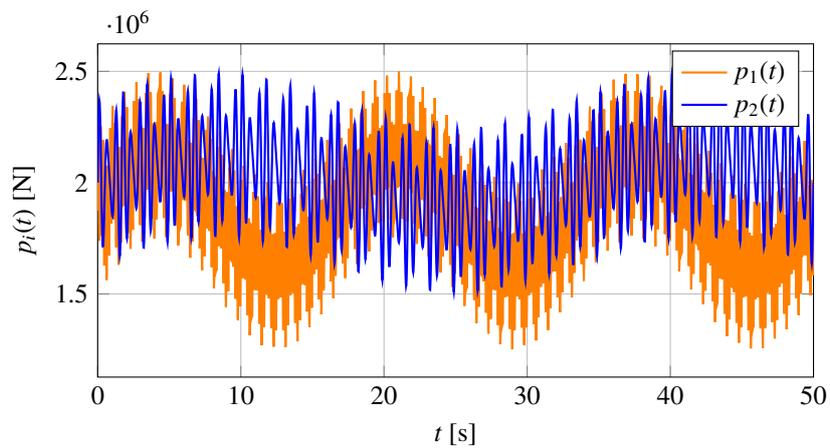
 
    \centering
    \include{elements/loadsCases_example3_plot}
    \setlength{\abovecaptionskip}{-0.3cm}
    \setlength{\belowcaptionskip}{-0.2cm}
    \caption{\small Numerical example of Sec. \ref{subsec:thirdex}. 50 seconds plot of the two load cases' magnitudes detailed in Table.~\ref{tab:loadCases_table_example3}.}
    \label{fig:loadCases_example3}
\end{figure*}

In this section, we discuss a structural optimization case characterized by a higher number of design variables and by stress constraints. The stress constraints considered limit the maximum tensile and compressive stresses in all bars, and the limit is defined by the yielding stress. Additional stress constraints are considered to limit stresses in compression, but this time the limit is defined by the stress associated to each bar Euler's buckling force. In the optimization process, the structural volume is minimized. We optimize the design of a truss structure, which is shown in Fig.~\ref{fig:systemDrawing_example3}. It is made of $13$ steel bars. We apply two separate sets of transient loads at the top free nodes of the structure, acting perpendicularly to the roof slope. On the left and right sides of the truss structure, the two sets of loads generate pressure and suction effects, respectively. While the time-history analysis is performed on the entire predefined time window, the optimization functions are evaluated only in a desired sub-time domain, which is specified later on. 
The information associated to the nodal coordinates of the structure and the elements' connectivity is provided in Table~\ref{tab:nodesData_example3} and Table~\ref{tab:elementsData_example3}.
\begin{table} [htbp!]
    \caption{\small Numerical example of Sec. \ref{subsec:thirdex}. Nodal coordinates of the truss structures shown in Fig.~\ref{fig:systemDrawing_example3}.}
    \centering
    \include{elements/nodesData_example3}
    \label{tab:nodesData_example3}
\end{table}
\begin{table} [htbp!]
    \caption{\small Numerical example of Sec. \ref{subsec:thirdex}. The elements' connectivity information for the truss structure considered (Fig.~\ref{fig:systemDrawing_example3}).}
    \centering
    \include{elements/elementsData_example3}
    \label{tab:elementsData_example3}
\end{table}
The two load cases considered acting on the external nodes of the structure are inspired by the dynamic loads used in \cite{verbart2018working}. They are detailed in Table~\ref{tab:loadCases_table_example3} and plotted in Fig.~\ref{fig:loadCases_example3}. 
\begin{table} 
    \caption{\small Numerical example of Sec. \ref{subsec:thirdex}. Values of the load cases considered, which are inspired by the loads used in \cite{verbart2018working}.}
    \centering
    \include{elements/loadCases_def_table_example3}
    \label{tab:loadCases_table_example3}
\end{table}
Note that to ensure zero displacements as initial conditions, the loads are multiplied by an exponential function in the form of $1-\exp(-0.2t)$, which increases the magnitude of the loads at the beginning of the structural response simulations starting from zero. This modification does not affect the structural behavior in the optimization time window considered to compute the system performance.

Given the axial symmetry of the structure, there are in total $7$ design variables, which are the diameters of structural elements' cross-sections. The link between the $7$ design variables and each $b$-th element of the structure, for $b=1,\dots, 13$, is defined in Eq.\ref{eq:variable2diam}: 
\begin{equation}
    \label{eq:variable2diam}
    D_b(\mathbf{x})=
    \begin{cases}
    x_1 &, b = 1,4 \\
    x_2 &, b = 2,3 \\
    x_3 &, b = 5,10 \\
    x_4 &, b = 6,9 \\
    x_5 &, b = 7,8 \\
    x_6 &, b = 11,13 \\
    x_7 &, b = 12
\end{cases}.
\end{equation}
In matrix form, this becomes:
\begin{equation}
    \label{eq:variable2diam_mat}
    \mathbf{D}(\mathbf{x})= \mathbf{L} \mathbf{x}
\end{equation}
where $L_{ij}=1$ if bar $i$ belongs to the size-group $j$, or in other words if the diameter $D_i$ equals $x_j$. Otherwise, $L_{ij}=0$. In this example, the matrix $\mathbf{L}$ has dimensions $[13 \times 7]$.
The Young's modulus of elasticity and the density of the bars are $E_b=200$ GPa and $\rho_b=7850$ kg/m$^3$, for $b=1,\dots13$. The elements are connected internally by hinges and the structure is anchored to the ground with full supports. Given that the loads are applied to the nodes, all bars experience only axial forces. We do not consider any structural damping.

Regarding the stress constraints, we consider a maximum allowable stress in tension and compression, defined as $\sigma_{max}=200$ MPa. The axial stress in each bar $b$ at time-step $k$ is:
\begin{equation}
    \label{eq:normalstress_calc}
    \sigma_{b,k}=E_b \varepsilon_{b,k} 
\end{equation}
where $\varepsilon$ is the strain in the bar, given by:
\begin{equation}
    \label{eq:strain_calc}
    \varepsilon_{b,k} = \frac{\Delta L_{b,k}}{L_{b,0}} = \frac{L_{b,k}-L_{b,0}}{L_{b,0}}.
\end{equation}
$L_{b,0}$ and $L_{b,k}$ are the initial and current lengths of the bar at time step $k$.  
The constraints are then:
\begin{equation}
    -\sigma_{max} \leq \sigma_{b,k} \leq \sigma_{max}.
\end{equation} 
The initial length of each bar is constant, but the length at time step $k$ depends on the bar nodal displacements. Denoting the $b$-th bar nodes as $n_{b1}$ and $n_{b2}$, and their displacements as $\mathbf{u}_k^{n_{b1}}$ and $\mathbf{u}_k^{n_{b2}}$, the locations of the nodes $\mathbf{X}_k^{n_{b1}}$ and $\mathbf{X}_k^{n_{b2}}$ can be computed as:
\begin{equation}
    \label{eq:nodesloc_calc}
    \begin{split}
    &\mathbf{X}_k^{n_{b1}} = \mathbf{X}_0^{n_{b1}} + \mathbf{u}_k^{n_{b1}} \\
    &\mathbf{X}_k^{n_{b2}} = \mathbf{X}_0^{n_{b2}} + \mathbf{u}_k^{n_{b2}}.
    \end{split}
\end{equation}
Finally, using Eq.~\eqref{eq:nodesloc_calc}, the bar's length can be computed and substituted into Eq.~\eqref{eq:strain_calc}:
\begin{equation}
    \label{eq:finallength_calc}
    L_{b,k} = \left\lVert \mathbf{X}_k^{n_{b1}}-\mathbf{X}_k^{n_{b2}} \right\rVert.
\end{equation}
In addition to the tension/compression stress constraints, the compression force in each bar is constrained to be less than the bars' buckling limit force. For a bar with a moment of inertia $I_{b,k}$ and length $L_{b}$, the compression force limit due to Euler buckling is:
\begin{equation}
    \label{eq:bucklingforce_calc}
    p_{b,k}^{buck} = \frac{\pi^2 E_b I_{b,k}}{L_{b}^2}
\end{equation}
with the moment of inertia given by:
\begin{equation}
    \label{eq:inertia_term_example2}
    I_b(\mathbf{x})=\frac{\pi D_b^4(\mathbf{x})}{64}.
\end{equation}
We thus consider additional stress constraints, by computing the equivalent buckling stress in each bar corresponding to the load of Eq.~\eqref{eq:bucklingforce_calc}:
\begin{equation}
    \label{eq:bucklingstress_calc}
    \sigma_{b,k}^{buck} = \frac{p_{b,k}^{buck}}{A_{b,k}} = \frac{\pi^2 E_b I_{b,k}}{L_b^2 A_{b,k}}.
\end{equation}
Note that the length $L_b$ considered in Eq.~\eqref{eq:bucklingstress_calc} is the initial length $L_{b,0}$. Writing more explicitly the dependency of the diameter of each bar to the actual optimization variables in Eq.~\eqref{eq:bucklingstress_calc}, we get the final expression for the buckling stress limit:
\begin{equation}
    \label{eq:bucklingstress_calc2}
    \sigma_{b,k}^{buck}(\mathbf{x}) = \frac{\pi^2 E_b}{16 L_b^2} D_b^2(\mathbf{x}).
\end{equation}
As a result, we consider also the following additional stress constraints in case of compression:
\begin{equation}
    -\sigma_{b,k}^{buck}(\mathbf{x}) \leq \sigma_{b,k}
\end{equation} 
where we assume that tension stresses, or forces, are positive and compression stresses are negative.

The reference dynamic response optimization problem formulation reads:
\begin{equation}
\label{eq:dynoptform_example3}
\begin{split}
    \underset{\mathbf{x} \in {\rm I\!R}^7}{\text{minimize }} &  f_c(\mathbf{x})=\sum_{b=1}^{13}{l_b \frac{\pi D_b^2(\mathbf{x})}{4}}\\
    \text{subject to:  } & g^l_{b,k}(\mathbf{x},\mathbf{u}^l_k(\mathbf{x}))= \frac{|\sigma_b(\mathbf{u}^l_k(\mathbf{x}))|}{\sigma_{max}} - 1 \leq 0 \\
    & q^l_{b,k}(\mathbf{x},\mathbf{u}^l_k(\mathbf{x}))= - 1-\frac{\sigma_b(\mathbf{u}^l_k(\mathbf{x}))}{\sigma_b^{buck}(\mathbf{x})} \leq 0 \\
    & x_i^{min} \leq x_i \leq x_i^{max}, \\ &i=1,\dots,7 \quad b=1,\dots,13 \\
    &k=s_t,\dots,n_t \quad l=1,2\\
    \text{with:  } &\mathbf{M} \ddot{\mathbf{u}}^l_h(\mathbf{x})
    +\mathbf{K}(\mathbf{x}) \mathbf{u}^l_h(\mathbf{x}) = \mathbf{p}_{h}^l\\
    & \text{for } h=1,2,\dots,n_t \quad l=1,2.
\end{split}
\tag{$\mathcal{P}^{3}_D$}
\end{equation}
In \eqref{eq:dynoptform_example3}, $l$ is the load case, $i$ and $b$ represent the design variable and bar indices, respectively. 
All structural members are modeled using Euler–Bernoulli beam finite elements \citep{przemieniecki1985theory}. The formulations for the consistent local mass and stiffness matrices and their derivatives with respect to a design variable are detailed in \ref{subsec:dynmat_ber} and \ref{subsec:diff_dynmat_ber}.
The time history analysis is performed for $h\in[0,n_t]$, but the constraints are evaluated only in the time window $k\in[s_t,n_t]$, with $t_{s_t}=30$ s and $t_{n_t}=50$ s. Within this discrete time window, there are $52052$ constraints that are considered. 

The F-ESL static response optimization problem formulation for the $W$-th iteration is:
\begin{equation}
\label{eq:statoptform_example3}
\begin{split}
    \underset{\mathbf{x} \in {\rm I\!R}^7}{\text{minimize }} &  f_c(\mathbf{x})=\sum_{b=1}^{13}{l_b \frac{\pi D_b^2(\mathbf{x})}{4}}\\
    \text{subject to:  } & g_{b,k}^{l}(\mathbf{x},\mathbf{\tilde{u}}_k^{l}(\mathbf{x}))= \frac{|\sigma_b(\mathbf{\tilde{u}}_k^{l}(\mathbf{x}))|}{\sigma_{max}} - 1 \leq 0 \\
    & q_{b,k}^{l}(\mathbf{x},\mathbf{\tilde{u}}_k^{l}(\mathbf{x}))= - 1-\frac{\sigma_b(\mathbf{\tilde{u}}_k^{l}(\mathbf{x}))}{\sigma_b^{buck}(\mathbf{x})} \leq 0 \\
    & x_i^{min} \leq x_i \leq x_i^{max}, \\ &i=1,\dots,7, \quad b=1,\dots,13 \\
    &k=s_t,\dots,n_t, \quad l=1,2 \\
    \text{with:  } &\mathbf{K} (\mathbf{x}) \mathbf{\tilde{u}}_k^{l}(\mathbf{x}, \mathbf{x}_W) 
     = \mathbf{f}^{FESL,l}_{k} (\mathbf{x}, \mathbf{x}_W)\\
    & \text{for } k=s_t,\dots,n_t \quad l=1,2 .
\end{split}
\tag{$\mathcal{P}^3_{S,W}$}
\end{equation}

The problem at hand was first solved directly by solving its original formulation defined in ~\eqref{eq:dynoptform_example3}, and then with the F-ESL approach by solving the sequence of static response sub-problems defined in ~\eqref{eq:statoptform_example3}. The F-ESL algorithm stopped once the stopping criteria of Eq.~\eqref{eq:stoppingcond} was satisfied, considering  $\epsilon=1\mathrm{e}{-6}$. Both \eqref{eq:dynoptform_example3} and \eqref{eq:statoptform_example3} were solved using the SQP algorithm of the \verb+fmincon+ function with its default settings, with the exception of \verb+ConstraintTolerance=1e-4+ and \verb+StepTolerance=1e-12+. The initial values for all design variables (i.e., bars' diameters) were set to $x_i^0=0.5$ m.

\begin{figure*}
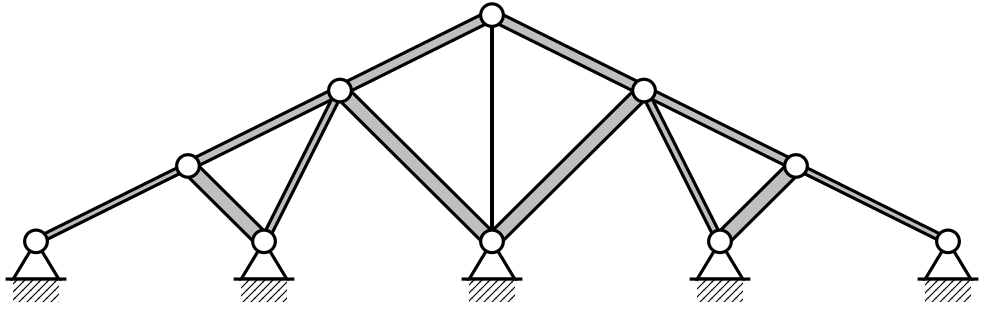
 
    \centering
    \include{elements/system_drawing_opt_example3}
    \setlength{\abovecaptionskip}{-0.3cm}
    \setlength{\belowcaptionskip}{-0.2cm}
    \caption{\small Numerical example of Sec. \ref{subsec:thirdex}. The optimized symmetric truss design obtained by solving the problem using the F-ESL formulation. Note that an identical solution was obtained solving directly the original problem for dynamic response optimization.}
    \label{fig:systemDrawing_opt_example3}
\end{figure*}

Table~\ref{tab:resultstable_example3_1} lists the number of time-history analyses, the objective values and the design variables of the final solutions obtained solving the problem at hand directly  and with the proposed F-ESL approach. We can observe that the final values of the objective function and the design variables of both solutions are identical\footnote{At least in a engineering perspective, i.e., with a very small difference ($4$ digits after the decimal point).}. Both formulations converged to the same solution, but the direct solution of the problem \eqref{eq:dynoptform_example3} required $35$ time-history analyses while F-ESL required only $4$ time-history analyses. In addition, Table~\ref{tab:resultstable_example3_2} lists the entries of the gradient of the Lagrangian at the final solutions. The values of the gradient entries are of the order of $10^{-5}$ or smaller, indicating that the solution is indeed a stationary point. 
The final optimal truss structure is illustrated in Fig.~\ref{fig:systemDrawing_opt_example3}. The reader should  note that in correspondence of the final solution, the diameter of bar $13$, i.e., the vertical bar in the middle, is equal to the minimum value allowed by the box constraints $D_{13}=x_{min}=0.001$. 
This means that the bar does not contribute to the truss load baring system and in principle could be removed. 
The axial stress of each bar is plotted for the first and second load cases in Fig.~\ref{fig:sigma_plot_loadcase1_example3} and Fig.~\ref{fig:sigma_plot_loadcase1_example3}, respectively, without considering bar $13$. On each plot, the stresses of the bars with the same design variable are plotted together. The stress limits $\pm \sigma_{max}$ and $\sigma_{buck}$ (buckling stress) are marked with red and black dashed lines, respectively. All the optimized stresses are within the prescribed limits. 

\begin{table*} [htbp!]\setlength{\tabcolsep}{4pt}
    \caption{\small Numerical example of Sec. \ref{subsec:thirdex}. Final solution results obtained with F-ESL, and by solving the original problem directly. The table provides the number of time-history-analyses performed $N_{THA}$, the final design $\mathbf{x}^*$, the final value of the objective function $f_c(\mathbf{x}^*)$.}
    \centering\scriptsize
    \include{elements/results_table_example3_1}
    \label{tab:resultstable_example3_1}
\end{table*}

\begin{table*} [htbp!]\setlength{\tabcolsep}{4pt}
    \caption{\small Numerical example of Sec. \ref{subsec:thirdex}. Final solution results obtained with F-ESL, and by solving the original problem directly. The table provides the  maximum stress constraint $g_{\max}=\max ([g_{b,k}^{l}(\mathbf{x}^*);q_{b,k}^{l}(\mathbf{x}^*)])$ and the gradient of the Lagrangian of the problem $\nabla \mathcal{L}(\mathbf{x}^*)$.}
    \centering\scriptsize
    \include{elements/results_table_example3_2}
    \label{tab:resultstable_example3_2}
\end{table*}

\begin{figure*} 
    \centering
    \include{elements/sigma_plot_loadcase1_example3}
    \setlength{\abovecaptionskip}{-1cm}
    \setlength{\belowcaptionskip}{-1cm}
    \caption{\small Numerical example of Sec. \ref{subsec:thirdex}. Axial stresses of the optimized structure in various bars for the first load case $p_1(t)$. The stress limits $\pm \sigma_{max}$ and $\sigma_{buck}$ (buckling stress) are marked with red and black dashed lines, respectively.}
    \label{fig:sigma_plot_loadcase1_example3}
\end{figure*}

\begin{figure*}
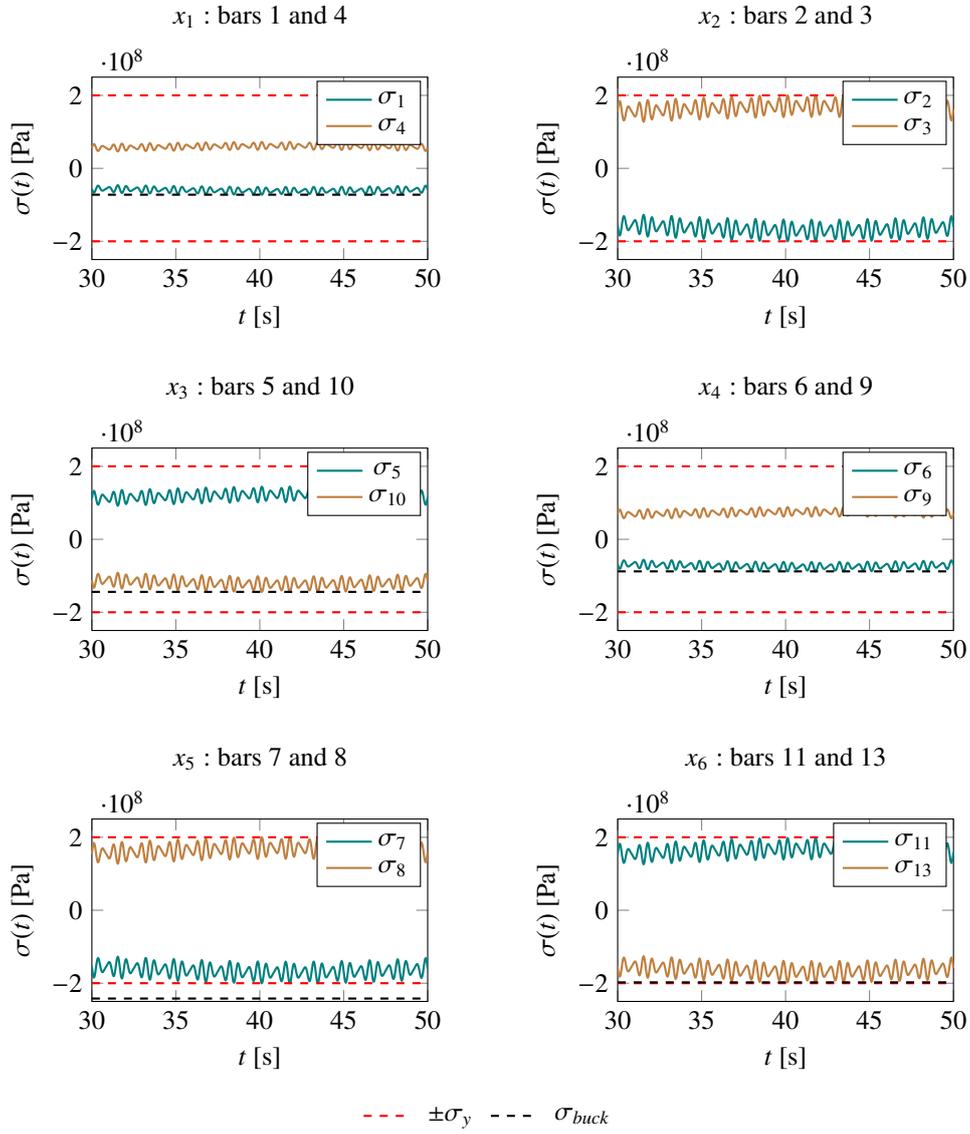
 
    \centering
    \include{elements/sigma_plot_loadcase2_example3}
    \setlength{\abovecaptionskip}{-0.3cm}
    \setlength{\belowcaptionskip}{-0.2cm}
    \caption{\small Numerical example of Sec. \ref{subsec:thirdex}. Axial stresses of the optimized structure in various bars for the second load case $p_2(t)$. The stress limits $\pm \sigma_{max}$ and $\sigma_{buck}$ (buckling stress) are marked with red and black dashed lines, respectively.}
    \label{fig:sigma_plot_loadcase2_example3}
\end{figure*}


\section{Conclusions}
\label{sec:end}
This paper introduces a novel first-order equivalent static loads (F-ESL) approach for structural optimization of dynamic response. F-ESL overcomes the known limitations of the original ESL formulation by including first-order information in the equivalent static loads definition. As a result, F-ESL is able to recognize final solutions that satisfy first-order optimality conditions of the original problem. 
By design,  F-ESL is designed to ensure consistency between the displacement sensitivities of the approximate static response and the original dynamic response. This consistency allows F-ESL to identify optimal solutions more accurately than traditional equivalent static loads formulations for structures subjected to dynamic loads. 

The proposed approach has been validated through several reproducible numerical examples. The first example considered consists of a simple one degree-of-freedom system with two variables. To further assess the proposed F-ESL approach, two additional, more complex, examples were studied. These are the design optimization of a shear frame subjected to a realistic earthquake, defined by the associated ground acceleration record, and the optimization of a truss structure with a larger number of design variables, two separate loading conditions, and stress and buckling constraints. Each example consistently demonstrated the capability of F-ESL to obtain the same optimized designs as those obtained solving directly the dynamic response optimization problems, and requiring less computational efforts (measured in terms of number of time-history analyses).

In all numerical applications studied, F-ESL converged to solutions closely matching those obtained directly solving the original dynamic response optimization problem, but requiring fewer transient analyses.  
In general, the results show that F-ESL is a reliable and effective approach for design optimization of structures subjected to transient loads. The results also show that F-ESL is capable of recognizing optimal designs of the original dynamic problem.

\appendix

\section{Sensitivity analysis}
\label{sec:appendixA}

\subsection{Direct sensitivity analysis}
\label{subsec:directsen_andx}
The Direct Sensitivity Analysis Method is a method based on the computation of the response sensitivities using direct differentiation of the governing equilibrium equations with respect to the design variables. In the case of structural dynamic response, we have the following dynamic equilibrium equations:
\begin{equation}
    \label{eq:dynamic_ODE}
    \mathbf{M} \ddot{\mathbf{u}}_k + \mathbf{C} \dot{\mathbf{u}}_k + \mathbf{K} {\mathbf{u}}_k - \mathbf{p}_k = \mathbf{0}.
\end{equation}
After differentiating both sides of Eq.~\eqref{eq:dynamic_ODE} by the design variable $x_i$ we obtain:
\begin{equation}
    \label{eq:dynamic_ODE_diff_xi}
    \mathbf{M} \frac{\partial \ddot{\mathbf{u}}_k}{\partial x_i}  + \mathbf{C} \frac{\partial \dot{\mathbf{u}}_k}{\partial x_i}  + \mathbf{K} \frac{\partial {\mathbf{u}}_k}{\partial x_i} = - \frac{\partial \mathbf{M}}{\partial x_i} \ddot{\mathbf{u}}_k - \frac{\partial \mathbf{C}}{\partial x_i} \dot{\mathbf{u}}_k - \frac{\partial \mathbf{K}}{\partial x_i} {\mathbf{u}}_k + \frac{\partial \mathbf{p}_k}{\partial x_i}.
\end{equation}
Assuming that the definition of the structural matrices is known, after solving the equations of motion of Eq.~\eqref{eq:dynamic_ODE}, using for example the Newmark-$\beta$ algorithm, the right-hand side is known and can be denoted as $\bar{\mathbf{p}}$:
\begin{equation}
    \label{eq:directsen_p_bar}
    \bar{\mathbf{p}} = - \frac{\partial \mathbf{M}}{\partial x_i} \ddot{\mathbf{u}}_k - \frac{\partial \mathbf{C}}{\partial x_i} \dot{\mathbf{u}}_k - \frac{\partial \mathbf{K}}{\partial x_i} {\mathbf{u}}_k + \frac{\partial \mathbf{p}_k}{\partial x_i}.
\end{equation}
Let us denote the displacement sensitivity of each time-step $k$ as $\mathbf{\lambda}_k=\frac{\partial \mathbf{u}_k}{\partial x_i}$, then Eq.~\eqref{eq:dynamic_ODE_diff_xi} turns into:
\begin{equation}
    \label{eq:directsen_ODE_lambda}
    \mathbf{M} \ddot{\bm{\lambda}}_k + \mathbf{C} \dot{\bm{\lambda}}_k + \mathbf{K} {\bm{\lambda}}_k = \bar{\mathbf{p}}_k.
\end{equation}
Since the initial conditions $\mathbf{u}_0$ and $\dot{\mathbf{u}}_0$ of Eq.~\eqref{eq:dynamic_ODE} are constant, the initial conditions in Eq.~\eqref{eq:directsen_ODE_lambda} are
\begin{equation}
    \label{eq:directsen_ODE_initialcond}
    \bm{\lambda}_0 = \bm{0}, \space \dot{\bm{\lambda}}_0 = \mathbf{0}.
\end{equation}
After solving Eq.~\eqref{eq:directsen_ODE_lambda} with a numerical integration solver similar to that used for Eq.~\eqref{eq:dynamic_ODE}, we obtain each $\bm{\lambda}_k$, which are the desired displacement sensitivities, $\frac{\partial \mathbf{u}_k}{\partial x_i}$. These are then used to compute the full gradients of the objective and constraints' functions.

\subsection{Gradient of the compliance objective function}
\label{subsec:compobjsen_andx}
The sensitivity for the compliance objective function in the first example has the same derivation for both the original dynamic response optimization problem ~\eqref{eq:dynoptform_example1} and the static response optimization sub-problems in the ESL and F-ESL approaches ~\eqref{eq:statoptform_example1}, up to the displacement sensitivities $\frac{\partial u_k}{\partial x_i}$. Differentiating the compliance in \eqref{eq:dynoptform_example1} for a specific time-step $k$ with respect to an arbitrary design variable $x_i$ leads to:
\begin{equation}
\label{eq:objsen_example1}
\begin{split}
 &\frac{\partial f_c^k}{\partial x_i} = \frac{\partial}{\partial x_i} \left(k u_k^2 \right) = \frac{\partial k}{\partial x_i} u_k^2 + 2ku_k\frac{\partial u_k}{\partial x_i},
\end{split}
\end{equation}
where the displacement sensitivities are computed using the direct sensitivity analysis method, detailed in \ref{subsec:directsen_andx}.

\subsection{Gradient of the drift constraints}
\label{subsec:driftconstsen_andx}
Following the derivations in Sec.~\ref{subsec:secondex}, differentiating the drift term in Eq.~\eqref{eq:driftcalc} by the design variable $x_i$ at time-step $k$ yields:
\begin{equation}
    \label{eq:driftsen}
    \frac{\partial \mathbf{d}_k}{\partial x_i} = \mathbf{D} \frac{\partial \mathbf{u}_k}{\partial x_i}
\end{equation}
since the drift matrix $\mathbf{D}$ is constant, i.e. $\partial\mathbf{D}/\partial x_i=\mathbf{0}$.
The displacement sensitivities $\frac{\partial\mathbf{u}_k}{\partial x_i}$ are computed using the Direct Sensitivity Analysis method, detailed in \ref{subsec:directsen_andx}, which requires the sensitivities of the dynamic matrices and of the external forces. The mass matrix in Eq.~\eqref{eq:stiff_mass_matrices} and the earthquake excitation force in Eq.~\eqref{eq:earthquakeforce} are constants, therefore:
\begin{equation}
    \label{eq:mass_earthquakeforce_diff}
    \frac{\partial\mathbf{M}}{\partial x_i}=\begin{bmatrix}
        0 & 0 \\
        0 & 0 \\
    \end{bmatrix}, \space \space \space \space 
    \frac{\partial\mathbf{f}_{g,k}}{\partial x_i}=\begin{bmatrix}
        0  \\
        0 \\
    \end{bmatrix}.
\end{equation}
The stiffness term in Eq.~\eqref{eq:columnStiffness} depends on the design variables through the area moment of inertia term in Eq.~\eqref{eq:inertia_term}, and therefore its derivative is:
\begin{equation}
    \label{eq:stiffness_diff}
    \frac{\partial k_i}{\partial x_i}=\frac{3 E_i}{L_i^3} \frac{\partial I_i(x_i)}{\partial x_i} = \frac{3 E_i}{L_i^3} \frac{\pi x_i^3}{16} 
\end{equation}
Note that $\frac{\partial k_i}{\partial x_j}=0$ for $i \neq j$. 
Substituting Eq.~\eqref{eq:stiffness_diff} in Eq.~\eqref{eq:stiff_mass_matrices} yields:
\begin{equation}
    \label{eq: stiff_matrix_diff}
    \frac{\partial \mathbf{K}}{dx_1}=\begin{bmatrix}
    2 \frac{\partial k_1}{\partial x_1}  & 0 \\
    0 & 0 \\
    \end{bmatrix},
    \quad
    \frac{\partial\mathbf{K}}{\partial x_2}=\begin{bmatrix}
    2 \frac{\partial k_2}{\partial x_2}  & -2 \frac{\partial k_2}{\partial x_2} \\
    -2 \frac{\partial k_2}{\partial x_2} & 2 \frac{\partial k_2}{\partial x_2}
    \end{bmatrix}.
\end{equation}
Differentiating the Rayleigh damping matrix in Eq.~\eqref{eq:damping_matrix}, we get the expression:
\begin{equation}
    \label{eq: damping_matrix_diff}
    \frac{d \mathbf{C}_s}{dx_i} = \frac{d \alpha_M}{dx_i} \mathbf{M} + \alpha_M \frac{d \mathbf{M}}{dx_i} + \frac{d \alpha_K}{dx_i} \mathbf{K} + \alpha_K \frac{d \mathbf{K}}{dx_i}
\end{equation}
Using the definitions of the natural frequencies in Eq.~\eqref{eq:eigenvalueprob} and Eq.~\eqref{eq:rayleighcoeff}, following the same procedure as in \cite{idels2020performance}, gives the sensitivities for the mass-proportional  Rayleigh damping coefficient:
\begin{equation}
    \label{eq:damping_coeff_aM_diff}
    \frac{\partial \alpha_M}{\partial x_i} = \xi \frac{2 \omega_2^2}{(\omega_1+\omega_2)^2} \frac{\partial \omega_1}{\partial x_i} + \xi \frac{2 \omega_1^2}{(\omega_1+\omega_2)^2} \frac{\partial \omega_2}{\partial x_i}
\end{equation}
and for the stiffness-proportional Rayleigh damping coefficient:
\begin{equation}
    \label{eq:damping_coeff_aK_diff}
    \frac{\partial \alpha_K}{\partial x_i} = -\xi \frac{2}{(\omega_1+\omega_2)^2} \frac{\partial \omega_1}{\partial x_i} - \xi \frac{2 }{(\omega_1+\omega_2)^2} \frac{\partial \omega_2}{\partial x_i}, 
\end{equation}
with:
\begin{equation}
    \label{eq:omega_diff}
    \frac{\partial \omega_j}{\partial x_i} = \frac{\mathbf{\hat{\phi}}_j^T (\frac{\partial \mathbf{K}}{\partial x_i}-\omega_j^2\frac{\partial \mathbf{M}}{\partial x_i}) \mathbf{\hat{\phi}}_j}{2 \omega_j}
\end{equation}
In Eq.~\eqref{eq:omega_diff}, $\mathbf{\hat{\phi}}_j$ are the normalized eigenvectors of Eq.~\eqref{eq:eigenvalueprob}.

\subsection{Gradient of the axial stress constraints}
\label{subsec:stressconstsen_andx}
Differentiating the axial stress term in Eq.~\eqref{eq:normalstress_calc} of a bar $b$ by the design variable $x_i$ at time-step $k$ gives:
\begin{equation}
    \label{eq:normailstress_diff}
    \frac{\partial \sigma_{b,k}}{\partial x_i} = E_b \frac{\partial \varepsilon_{b,k}}{\partial x_i} = \frac{E_b}{L_{b,0}} \frac{\partial L_{b,k}}{\partial x_i}.
\end{equation}
Using Eq.~\eqref{eq:nodesloc_calc} and Eq.~\eqref{eq:finallength_calc}, we can differentiate the length of the bar as follows:
\begin{equation}
    \label{eq:finallength_diff}
    \frac{\partial L_{b,k}}{\partial x_i} = \frac{(\mathbf{X}_k^{n_{b1}}-\mathbf{X}_k^{n_{b2}})^T (\frac{\partial\mathbf{u}_k^{n_{b1}}}{\partial x_i}-\frac{\partial\mathbf{u}_k^{n_{b2}}}{\partial x_i})}{\left\lVert \mathbf{X}_k^{n_{b1}}-\mathbf{X}_k^{n_{b2}} \right\rVert}.
\end{equation}
Substituting Eq.~\eqref{eq:finallength_diff} into Eq.~\eqref{eq:normailstress_diff}, we get the final expression for the normal stress sensitivities, where the displacements sensitivities $\frac{\partial\mathbf{u}_k^{n_{b1}}}{\partial x_i}$ and $\frac{\partial\mathbf{u}_k^{n_{b2}}}{\partial x_i}$ are computed separately, as described in \ref{subsec:directsen_andx}.

\subsection{Gradient of the buckling stress constraints}
\label{subsec:buckstressconstsen_andx}
Differentiating the buckling stress term in Eq.~\eqref{eq:bucklingstress_calc2} of bar $b$ by the design variable $x_i$ at time-step $k$ gives:
\begin{equation}
    \label{eq:bucklingstress_diff}
    \frac{\partial \sigma_{b,k}^{buck}(\mathbf{x})}{\partial x_i} = \frac{\pi^2 E_b}{8 L_b^2} D_b(\mathbf{x}).
\end{equation}
Once the buckling stress sensitivities of Eq.~\eqref{eq:bucklingstress_diff} are computed, the sensitivities of the buckling stress constraints considered in problems \eqref{eq:dynoptform_example3} and \eqref{eq:statoptform_example3} are calculated as follows:
\begin{equation}
    \label{eq:bucklingstress_cond_diff}
    \frac{\partial q_{b,k}}{\partial x_i} = \frac{\sigma_b \frac{\partial \sigma_b^{buck}}{\partial x_i} - \frac{\partial\sigma_b}{\partial x_i} \sigma_b^{buck}}{(\sigma_b^{buck})^2}.
\end{equation}

\section{Mass and stiffness matrices of a Euler-Bernoulli beam element}
\label{sec:appendixB}

\subsection{Mass and stiffness matrices}
\label{subsec:dynmat_ber}
The consistent mass and stiffness matrices used in Sec.~\ref{subsec:thirdex} are those of the Euler-Bernoulli beam elements \citep{przemieniecki1985theory}. In a 2D structure, like the ones considered here, the matrices presented have dimensions $[6 \times 6]$, due to the presence of three degrees of freedom for each node of the beam. 
For an element $e$, the element matrices depend on the design variables (beams' diameters).
Thus, the element mass matrix is: 
\begin{equation}
    \label{eq:mass_mat_example3}
    \begin{gathered}
    \mathbf{M}_e(x) = \\
    \frac{\rho A(x) l}{420}
    \begin{bmatrix}
    140 &   0 &   0 &  70 &   0 &    0 \\
      0 & 156 &  22l &   0 &  54 & -13l \\
      0 & 22l &  4l^2 &   0 & 13l & -3l^2 \\
     70 &   0 &   0 & 140 &   0 &    0 \\
      0 &  54 & 13l &   0 & 156 & -22l \\
      0 & -13l & -3l^2 &  0 & -22l & 4l^2
    \end{bmatrix}.
    \end{gathered}
\end{equation}
The element stiffness matrix is: 
\begin{equation}
    \label{eq:stiff_mat_example3}
    \begin{gathered}
    \mathbf{K}_e(x) = \\
    \frac{E}{8l^3}
    {\setlength{\arraycolsep}{3pt}\begin{bmatrix}
    8A l^2 & 0 & 0 & -8A l^2 & 0 & 0 \\
    0 & 96 I & 48 I l & 0 & -96 I & 48 I l \\
    0 & 48 I l & 32 I l^2 & 0 & -48 I l & 16 I l^2 \\
    -8A l^2 & 0 & 0 & 8A l^2 & 0 & 0 \\
    0 & -96 I & -48 I l & 0 & 96 I & -48 I l \\
    0 & 48 I l & 16 I l^2 & 0 & -48 I l & 32 I l^2
    \end{bmatrix}}
    \end{gathered},
\end{equation}
where $l$ is the length of the beam, $\rho$ is the density, $E$ is the modulus of elasticity. In Eq.~\eqref{eq:stiff_mat_example3} $A$ and $I$ depend on $x$: $A(x)=\pi D(x)^2/4$ is the cross-section area, $I(x)=\pi D(x)^4/64$ is the moment of inertia of the beam in the $z$ direction (perpendicular to the plane).

\subsection{Gradient of the mass and stiffness matrices}
\label{subsec:diff_dynmat_ber}
Direct differentiation of the mass and stiffness matrices with respect to a design variable $x$ gives:
\begin{equation}
    \label{eq:diff_mass_mat_example3}
    \begin{gathered}
    \frac{\partial \mathbf{M}_e(x)}{\partial x} = \\
    \frac{\pi \rho l x}{420}
    \begin{bmatrix}
    70 & 0 & 0 & 35 & 0 & 0 \\
    0 & 78 & 11l & 0 & 27 & -13l \\
    0 & 11l & 2l^2 & 0 & 13l & -3l^2 \\
    35 & 0 & 0 & 70 & 0 & 0 \\
    0 & 27 & 13l & 0 & 78 & -11l \\
    0 & -13l & -3l^2 & 0 & -11l & 2l^2
    \end{bmatrix}
    \end{gathered}
\end{equation}
and 
\begin{equation}
    \label{eq:diff_stiff_mat_example3}
    \begin{gathered}
    \frac{\partial \mathbf{K}_e(x)}{\partial x} = \\
    \frac{\pi E}{8L^3}
    \begin{bmatrix}
    4x l^2 & 0 & 0 & -4x l^2 & 0 & 0 \\
    0 & 6x^3 & 3x^3 l & 0 & -6x^3 & 3x^3 l \\
    0 & 3x^3 l & 2x^3 l^2 & 0 & -3x^3 l & x^3 l^2 \\
    -4x l^2 & 0 & 0 & 4x l^2 & 0 & 0 \\
    0 & -6x^3 & -3x^3 l & 0 & 6x^3 & -3x^3 l \\
    0 & 3x^3 l & x^3 l^2 & 0 & -3x^3 l & 2x^3 l^2
    \end{bmatrix}.
    \end{gathered}
\end{equation}
These matrices are used to calculate the sensitivities of the stress constraints, as shown in detail in \ref{subsec:stressconstsen_andx} and \ref{subsec:buckstressconstsen_andx}.

\section*{ Replication of results}
\noindent The numerical and algorithmic details are provided
in sufficient detail in the text to be implemented. 

\section*{ Conflict of interest}
\noindent The authors declared that they have no conflict of interest.

\section*{ Acknowledgments}
\noindent NP acknowledges the support of the Israel Science Foundation (ISF), grant No. 304/24.


\bibliography{myBib} 

\end{document}

%% file: elements/pseudo_algo.tex
{
\setlength{\baselineskip}{0.8cm}
\newcommand{\vspacealgo}{\vspace{0.8em}}

\begin{algorithmic}
\State Initialize  $\mathbf{x}_0\in \mathcal{X}$ 
\State $W \gets 1$
\State $\mathbf{x}_W \gets \mathbf{x}_0$
    \vspacealgo
\While{$||\mathbf{x}_{W}-\mathbf{x}_{W-1}||\leq \epsilon$} 
    \vspacealgo
    \State Perform dynamic-response analysis 
    \Statex \quad \quad $\forall k \in [1,N_t]$ with $\mathbf{x}=\mathbf{x}_W$
    \vspacealgo
    \State Compute equivalent static loads:
    \State \quad If $ESL$: use Eq.~\eqref{eq:esldef} 
    \State \quad If $F\text{-}ESL$: use Eq.~\eqref{eq:feslfdef}
    \State Solve the linear static problem \eqref{eq:statoptform}
    \State Denote the solution found $\mathbf{\tilde{x}}_W$
        \vspacealgo
    \State Variables update $\mathbf{x}_{W+1} \gets \mathbf{\tilde{x}}_W$
    \State Iteration update $W \gets W+1$
    \EndWhile
\end{algorithmic}
}



%% file: elements/system_drawing_example1.tex
\begin{tikzpicture}[every node/.style={draw,outer sep=0pt,thick}]
    \tikzstyle{ground_ver}=[fill, pattern=north east lines, draw=none,minimum width=1cm, minimum height=0.6cm]
    
    \tikzstyle{ground_hor}=[fill, pattern=north east lines, draw=none,minimum width=0.6cm, minimum height=2cm]
    
    \def \beamHeight {3cm}
    \def \beamWidth {3.5mm}

    \node (beam1) [fill=lightgray, very thick, minimum width=\beamWidth,minimum height=\beamHeight] {};
    \node (beam2) at (beam1.south) [fill=lightgray, very thick, minimum width=\beamWidth,minimum height=(1.5*\beamHeight), yshift=-(1.5*\beamHeight)/2] {};
    
    \node (ground_down) at (beam2.south) [ground_ver, yshift=0cm, xshift=0cm, anchor=north] {};
    \draw [very thick] (ground_down.north west) -- (ground_down.north east);
    
    \node (ground_up) at (beam1.north) [ground_ver, yshift=0cm, xshift=0cm, anchor=south] {};
    \draw [very thick] (ground_up.south west) -- (ground_up.south east);
    
    \node (ground_middle) at (beam1.south) [ground_hor, yshift=0cm, xshift=-1cm, anchor=east] {};
    \draw [very thick] (ground_middle.south east) -- (ground_middle.north east);
    
    \draw [very thick] (ground_middle.east)  + (0.3cm+1/3cm,-0.35) -- +(0.3cm+1/3cm,0.35cm) -- + (-0.15cm+1cm,0) -- + (0.3cm+1/3cm,-0.35);
    
    \draw [very thick] (ground_middle.east) ++ (0.15cm,-0.25cm) circle (0.15cm)  (ground_middle.east) ++ (0.15cm, 0.25cm) circle (0.15cm);
    
    \draw [fill=white, very thick] (beam1.south) ++ (0cm, 0cm) circle (0.2cm);
    \draw [-latex, ultra thick] (beam1.south) ++(0.2cm, 0cm) -- +(\beamWidth, 0cm) -- +(\beamWidth,-1/3*1.5*\beamHeight);
    
    \node [anchor=west, draw=none] (beam1text) at ($(beam1.east)+(0.15cm,0cm)$) {$\rho_1, E_1, l_1$};
    
    \node [anchor=west, draw=none] (beam2text) at ($(beam2.east)+(0.15cm,0cm)$) {$\rho_2, E_2, l_2$};
    
    \node [anchor=west, draw=none] (forcetext) at ($(beam1.south east)+(0.5cm,-1/3*1.5*1/2*\beamHeight)$) {$r(t) = r_0 \sin(\omega t)$};

    \node[anchor=east, draw, circle, fill=white, inner sep=1pt] (beam1number) at ($(beam1.west)-(0.2cm,0cm)$) {1};

    \node[anchor=east, draw, circle, fill=white, inner sep=1pt] (beam1number) at ($(beam2.west)-(0.2cm,0cm)$) {2};
\end{tikzpicture}

%% file: elements/results_table_example1.tex
\begin{tabular}{lcccccccc} \toprule\normalsize
 & $N_{THA}$&    $\mathbf{x}^{*}$&  $f_c(\mathbf{x}^{*})$& $\frac{\partial f_c(\mathbf{x}^{*})}{\partial x_1}$& $\frac{\partial f_c(\mathbf{x}^{*})}{\partial x_2}$& $g(\mathbf{x}^*)$& $\frac{\partial \mathcal{L}(\mathbf{x}^{*})}{\partial x_1}$& $\frac{\partial \mathcal{L}(\mathbf{x}^{*})}{\partial x_2}$ 
\\ \midrule
\eqref{eq:dynoptform_example1}& $7$&  $(0.1,0.6)$&  $0.888$&  $1.26$& $-1.56$& $-1.1\mathrm{e}{-16}$& $9.6\mathrm{e}{-13}$& $2.8\mathrm{e}{-12}$ 
\\
\eqref{eq:statoptform_example1}$_{ESL}$&  $2$& $(0.85,0.1)$&  $7.344$&  $-0.72$&  $-38.86$& $-1.1\mathrm{e}{-16}$& $7.29$& $-33.52$
\\
\eqref{eq:statoptform_example1}$_{F\text{-}ESL}$& $3$&  $(0.1,0.6)$&  $0.888$&   $1.26$& $-1.56$& $-1.1\mathrm{e}{-16}$& $-4.4\mathrm{e}{-16}$& $-8.9\mathrm{e}{-16}$ \\ \bottomrule
\end{tabular}

%% file: elements/VariablePath_plot_example1.tex
\begin{tikzpicture}
    \begin{axis}[
        colormap/blackwhite,
        xlabel={$x_1$}, ylabel={$x_2$},
        width=10cm,
        height=10cm,
        xmin = 0.05, xmax = 1.05,  
        ymin = 0.05, ymax = 1.05, 
    ]
         \addplot[dotted, very thick, forget plot] coordinates {
            (0.1, 0.1)
            (1.0, 0.1)
            (1.0, 1.0)
            (0.1, 1.0)
            (0.1, 0.1) 
        };

         \addplot[dotted, very thick, draw=red, forget plot] table[x index=0, y index=1] {data/VolumeConstData.dat};
        
	\addplot[contour prepared, draw=black, thick, forget plot] file {data/contourData.dat};

        \addplot [
                fill=gray,
                fill opacity=0.2, 
                draw=none, 
                forget plot
            ] table [x index=0, y index=1] {data/FeasibleSetData.dat};

        \addplot[very thick, black, solid,
             mark=*, mark options={fill=black}] table[x index=0, y index=1] {data/pathData_dyn.dat};
        \addplot[
            very thick,
            black,
            solid,
            postaction={
                decorate,
                decoration={
                    markings,
                    mark=at position 0.65 with {\arrow[scale=1.2]{stealth}}, 
                }
            }, 
            forget plot
        ] table[x index=0, y index=1] {data/pathData_dyn.dat};
        \addlegendentry{Reference}

        \addplot[very thick, teal,
             mark=*, mark options={fill=teal}] table[x index=0, y index=1] {data/pathData_ESL.dat};
        \addplot[
            very thick,
            teal,
            postaction={
                decorate,
                decoration={
                    markings,
                    mark=at position 0.5 with {\arrow[scale=1.2]{stealth}},
                }
            }, 
            forget plot
        ] table[x index=0, y index=1] {data/pathData_ESL.dat};
        \addlegendentry{$ESL$}

        \addplot[very thick, brown,
             mark=*, mark options={fill=brown}] table[x index=0, y index=1] {data/pathData_fesl.dat};
        \addplot[
            very thick,
            brown,
            postaction={
                decorate,
                decoration={
                    markings,
                    mark=at position 0.4 with {\arrow[scale=1.5]{stealth}},
                    mark=at position 0.9 with {\arrow[scale=1.2]{stealth}}
                }
            }, 
            forget plot
        ] table[x index=0, y index=1] {data/pathData_fesl.dat}; 
        \addlegendentry{$F\text{-}ESL$}

        \addplot[only marks, mark=*, black] coordinates {(0.1, 0.6)};
        \node[anchor=south west] at (axis cs:0.1, 0.6) {$\mathbf{x}^{*}_{FESL}=\mathbf{x}^{*}_{opt}$};

        \addplot[only marks, mark=*, black] coordinates {(0.85, 0.1)};
        \node[anchor=south west] at (axis cs:0.85, 0.1) {$\mathbf{x}^{*}_{ESL}$};
        
	\end{axis}
\end{tikzpicture}

%% file: elements/dUdx_xfesl_plot.tex
  \begin{subfigure}[b]{0.9\textwidth}
    \centering
    \begin{tikzpicture}
    \begin{axis}[
        xlabel={$t$},
        ylabel={$\frac{\partial u(\mathbf{x}^*_{opt})}{\partial x_1}$},
        width=8cm,
        height=6cm,
        xmin=200, xmax=202,
        legend to name=sharedlegend,
        legend style={draw=none},
        grid
    ]
    \addplot[color=black, very thick] table[x index=0, y index=1] {data/dynamic1_xfesl_Data.dat};
    \addplot[color=teal, dashed, very thick] table[x index=0, y index=1] {data/esl1_xfesl_Data.dat};
    \addplot[color=brown, dashed, very thick] table[x index=0, y index=1] {data/fesl1_xfesl_Data.dat};

    \addlegendentry{Dynamic}
    \addlegendentry{$ESLM$}
    \addlegendentry{$F\text{-}ESLM$}
    \end{axis}
    \end{tikzpicture}
\end{subfigure}

\vspace{1em}

\begin{subfigure}[b]{0.9\textwidth}
    \centering
    \begin{tikzpicture}
    \begin{axis}[
        xlabel={$t$},
        ylabel={$\frac{\partial u(\mathbf{x}^*_{opt})}{\partial x_2}$},
        width=8cm,
        height=6cm,
        xmin=200, xmax=202,
        grid
    ]
    \addplot[color=black, very thick] table[x index=0, y index=1] {data/dynamic2_xfesl_Data.dat};
    \addplot[color=teal, dashed, very thick] table[x index=0, y index=1] {data/esl2_xfesl_Data.dat};
    \addplot[color=brown, dashed, very thick] table[x index=0, y index=1] {data/fesl2_xfesl_Data.dat};
    \end{axis}
    \end{tikzpicture}
\end{subfigure}

    \begin{tikzpicture}
        \begin{axis}[
            hide axis,
            xmin=0,
            xmax=1,
            ymin=0,
            ymax=1,
            legend columns=3,
            legend style={draw=none, column sep=1ex},
        ]
        \addlegendimage{black, very thick}
        \addlegendentry{Reference}
        \addlegendimage{teal, dashed, very thick}
        \addlegendentry{$ESL$}
        \addlegendimage{brown, dashed, very thick}
        \addlegendentry{$F\text{-}ESL$}
        \end{axis}
    \end{tikzpicture}

%% file: elements/dUdx_xesl_plot.tex
    \begin{subfigure}[b]{0.9\textwidth}
        \centering
        \begin{tikzpicture}
        \begin{axis}[
            xlabel={$t$},
            ylabel={$\frac{\partial u(\mathbf{x}^*_{ESL})}{\partial x_1}$},
            width=8cm,
            height=6cm,
            xmin=200, xmax=202, 
            grid
        ]
        \addplot[color=black, very thick] table[x index=0, y index=1] {data/dynamic1_xesl_Data.dat};

        \addplot[color=teal, dashed, very thick] table[x index=0, y index=1] {data/esl1_xesl_Data.dat};

        \addplot[color=brown, dashed, very thick] table[x index=0, y index=1] {data/fesl1_xesl_Data.dat};
        \end{axis}
        \end{tikzpicture}
    \end{subfigure}
    
    \vspace{1em}
    
    \begin{subfigure}[b]{0.9\textwidth}
        \centering
        \begin{tikzpicture}
        \begin{axis}[
            xlabel={$t$},
            ylabel={$\frac{\partial u(\mathbf{x}^*_{ESL})}{\partial x_2}$},
            width=8cm,
            height=6cm,
            xmin=200, xmax=202, 
            transpose legend,
            legend columns=3,
            legend style={at={(0.5,-0.25)},anchor=north},
            grid
        ]
        \addplot[color=black, very thick] table[x index=0, y index=1] {data/dynamic2_xesl_Data.dat};

        \addplot[color=teal, dashed, very thick] table[x index=0, y index=1] {data/esl2_xesl_Data.dat};

        \addplot[color=brown, dashed, very thick] table[x index=0, y index=1] {data/fesl2_xesl_Data.dat};

        \end{axis}
        \end{tikzpicture}
    \end{subfigure}

    \begin{tikzpicture}
        \begin{axis}[
            hide axis,
            xmin=0,
            xmax=1,
            ymin=0,
            ymax=1,
            legend columns=3,
            legend style={draw=none, column sep=1ex},
        ]
        \addlegendimage{black, very thick}
        \addlegendentry{Reference}
        \addlegendimage{teal, dashed, very thick}
        \addlegendentry{$ESL$}
        \addlegendimage{brown, dashed, very thick}
        \addlegendentry{$F\text{-}ESL$}
        \end{axis}
    \end{tikzpicture}

    

%% file: elements/dFesl_dA.tex
\begin{tikzpicture}
    \begin{axis}[
        xlabel={$t$},
        ylabel={$\frac{\partial \mathbf{f}^{eq} (\mathbf{x}^*_{opt})}{\partial x}$},
        width=8cm,
        height=6cm,
        xmin=200, xmax=202, 
        xtick={200, 200.5, 201, 201.5, 202},
        grid,
        legend style={draw=none, at={(0.5,-0.3)}, anchor=north, legend columns=2}, 
    ]
    
    \addplot[color=blue, very thick] table[x index=0, y index=1] {data/dFesl_da1.dat};
    \addlegendentry{$\frac{\partial \mathbf{f}^{eq}(\mathbf{x}^*_{opt})}{\partial x_1}$}

    \addplot[color=red, very thick] table[x index=0, y index=1] {data/dFesl_da2.dat};
    \addlegendentry{$\frac{\partial \mathbf{f}^{eq}(\mathbf{x}^*_{opt})}{\partial x_2}$}

    \end{axis}
\end{tikzpicture}

%% file: elements/system_drawing_example2.tex
\begin{tikzpicture}[every node/.style={draw,outer sep=0pt,thick}]

    \tikzstyle{ground}=[fill, pattern=north east lines, draw=none, minimum width=0.5cm, minimum height=0.3cm]

    \def \storyHeight {2.5cm}
    \def \beamWidth {2mm}
    \def \frameWidth {2.5cm}

    \node (col1_bottom) [fill=lightgray, minimum width=\beamWidth, minimum height=\storyHeight, very thick] {};
    \node (col1_top) at (col1_bottom.north) [fill=lightgray, minimum width=\beamWidth, minimum height=\storyHeight, anchor=south, very thick] {};

    \node (col2_bottom) at ($(col1_bottom.east)+( \frameWidth, 0 )$) [fill=lightgray, minimum width=\beamWidth, minimum height=\storyHeight, very thick] {};
    \node (col2_top) at (col2_bottom.north) [fill=lightgray, minimum width=\beamWidth, minimum height=\storyHeight, anchor=south, very thick] {};

    \draw[very thick, fill=white] 
        ($(col1_bottom.north)+(-\beamWidth,0)$) -- 
        ($(col2_bottom.north)+(+\beamWidth,0)$) -- 
        ($(col2_bottom.north)-(-\beamWidth,2*\beamWidth)$) -- 
        ($(col1_bottom.north)-(+\beamWidth,2*\beamWidth)$) -- cycle;

    \draw[very thick, fill=white] 
        ($(col1_top.north)+(-\beamWidth,0)$) -- 
        ($(col2_top.north)+(+\beamWidth,0)$) -- 
        ($(col2_top.north)-(-\beamWidth,2*\beamWidth)$) -- 
        ($(col1_top.north)-(+\beamWidth,2*\beamWidth)$) -- cycle;

    \node (groundL) at (col1_bottom.south) [ground, anchor=north] {};
    \draw[very thick] (groundL.north west) -- (groundL.north east);

    \node (groundR) at (col2_bottom.south) [ground, anchor=north] {};
    \draw[very thick] (groundR.north west) -- (groundR.north east);

    \node[draw=none, anchor=center] at ($1/2*(col2_bottom.north)+1/2*(col1_bottom.north)+(0,2*\beamWidth)$) {$m_1$};
    \node[draw=none, anchor=center] at ($1/2*(col2_top.north)+1/2*(col1_top.north)+(0,2*\beamWidth)$) {$m_2$};

    \draw[-latex, ultra thick] ($(col2_bottom.north west) + (0.35cm, -\beamWidth)$) -- ($(col2_bottom.north west) + (1cm, -\beamWidth)$);
    \node[draw=none, anchor=west] at ($(col2_bottom.north west)+(0.2cm,\beamWidth)$) {$u_1$};

    \draw[-latex, ultra thick] ($(col2_top.north west) + (0.35cm, -\beamWidth)$) -- ($(col2_top.north west) + (1cm, -\beamWidth)$);
    \node[draw=none, anchor=west] at ($(col2_top.north west)+(0.2cm,\beamWidth)$) {$u_2$};    
    
    \node[draw=none, anchor=west] at ($(col2_bottom.east)$) {$l_1,E_1$};
    \node[draw=none, anchor=west] at ($(col2_top.east)$) {$l_2,E_2$};
    
    \node[anchor=east, draw, circle, fill=white, inner sep=1pt] at ($(col1_bottom.west)-(0.2cm,0)$) {1};
    \node[anchor=east, draw, circle, fill=white, inner sep=1pt] at ($(col1_top.west)-(0.2cm,0)$) {2};

    \node (groundAcc) at ($0.5*(col1_bottom.south) + 0.5*(col2_bottom.south)$) [anchor=center, draw=none] {};
    
    \node [anchor=center, draw=none] at (groundAcc) {$a_g$};
    
    \draw[latex-, ultra thick] ($(groundAcc) - (0.65cm, 0)$) -- ($(groundAcc) - (0.25cm, 0)$);
    
    \draw[-latex, ultra thick] ($(groundAcc) + (0.25cm, 0)$) -- ($(groundAcc) + (0.65cm, 0)$);

\end{tikzpicture}

%% file: elements/earthquake_plot.tex


\begin{tikzpicture}
    \begin{axis}[
        xlabel={$t$ [s]},
        ylabel={$a_g(t)$ [m/s$^2$]},
        width=10cm,
        height=6cm,
        xmin=0, xmax=20, 
        grid
    ]
    \addplot[color=blue, very thick] table[x index=0, y index=1] {data/agData.dat};

    \end{axis}
\end{tikzpicture}

%% file: elements/results_table_example2.tex
\begin{tabular}{llccccccccc} \toprule\small
 &  & $N_{THA}$ & $x^{*}_1$ & $x^{*}_2$ & $f_c(\mathbf{x}^{*})$ 
 & $\frac{\partial f_c(\mathbf{x}^{*})}{\partial x_1}$ & $\frac{\partial f_c(\mathbf{x}^{*})}{\partial x_2}$ 
 & $g_{\max}$ 
 & $\frac{\partial \mathcal{L}(\mathbf{x}^{*})}{\partial x_1}$ 
 & $\frac{\partial \mathcal{L}(\mathbf{x}^{*})}{\partial x_2}$ 
\\ \midrule
\eqref{eq:dynoptform_example2} & GA      & $10173$ & $0.2745$ & $0.2501$ & $0.7580$ & $3.02$ & $2.75$ & $-3.2\mathrm{e}{-5}$ & $-$ & $-$ \\
\eqref{eq:dynoptform_example2} & SQP   & $38$    & $0.2744$ & $0.2501$ & $0.7579$ & $3.02$ & $2.75$ & $-3.9\mathrm{e}{-15}$ & $4.6\mathrm{e}{-10}$ & $-6.8\mathrm{e}{-10}$ \\
\eqref{eq:statoptform_example2} & F-ESL & $15$    & $0.2744$ & $0.2501$ & $0.7579$ & $3.02$ & $2.75$ & $-1.5\mathrm{e}{-8}$ & $4.3\mathrm{e}{-6}$ & $3.1\mathrm{e}{-6}$ \\
\bottomrule
\end{tabular}

%% file: elements/drifts_plot_example2.tex

    
    


    


\begin{tikzpicture}
    \begin{axis}[
        xlabel={$t$ [s]},
        ylabel={$d_i(t)$ [m]},
        width=10cm,
        height=6cm,
        xmin=0, xmax=20, 
        ymin=-0.15, ymax=0.15, 
        grid
    ]
    
    \addplot[color=red, very thick, dotted, forget plot] coordinates {(0,0.1) (20,0.1)};
    
    \addplot[color=red, very thick, dotted,forget plot] coordinates {(0,-0.1) (20,-0.1)};

    \addplot[color=blue, very thick] table[x index=0, y index=1] {data/drifts_u1_example2.dat};
    \addlegendentry{$d_1(t)$}    

    \addplot[color=orange, very thick] table[x index=0, y index=1] {data/drifts_u2_example2.dat};
    \addlegendentry{$d_2(t)$}    
    
    \end{axis}
\end{tikzpicture}

%% file: elements/obj_iterations_plot_example2.tex

    


    


\begin{tikzpicture}
    \begin{axis}[
        xlabel={iteration},
        ylabel={$f_{c}(\mathbf{x}_W)$ [m$^3$]},
        width=10cm,
        height=6cm,
        xmin=0, xmax=40, 
        ymin=0, ymax=3, 
        grid
    ]
    
    \addplot[color=red, very thick, dotted, forget plot] table[x index=0, y index=1] {data/iterationsObj_conv_example2.dat};

    \addplot[color=black, very thick] table[x index=0, y index=1] {data/iterationsObj_example2.dat};
    \addlegendentry{Reference}

    \addplot[color=brown, very thick] table[x index=0, y index=1] {data/iterationsObj_fesl_example2.dat};
    \addlegendentry{$F\text{-}ESL$} 
    
    \end{axis}
\end{tikzpicture}

%% file: elements/VariablePath_plot_example2.tex
\begin{tikzpicture}
    \begin{axis}[
        colormap/blackwhite,
        xlabel={$x_1$}, ylabel={$x_2$},
        width=10cm,
        height=10cm,
        xmin = 0.15, xmax = 0.52,  
        ymin = 0.15, ymax = 0.52, 
        legend style={at={(0.02,0.98)}, anchor=north west}
    ]
         \addplot[dotted, very thick, forget plot] coordinates {
            (0.5, 0.15)
            (0.5, 0.5)
            (0.15, 0.5)
        };

         \addplot[dotted, very thick, draw=red, forget plot] table[x index=0, y index=1] {data/ConBoundaryData_example2.dat};
        
	\addplot[contour prepared, draw=black, thick, forget plot] file {data/contourData_example2.dat};

        \addplot [
                fill=gray,
                fill opacity=0.2, 
                draw=none, 
                forget plot
            ] table [x index=0, y index=1] {data/FeasibleSetData_example2.dat};

        \addplot[very thick, black, solid,
             mark=*, mark options={fill=black}] table[x index=0, y index=1] {data/pathData_dyn_example2.dat};
        \addplot[
            very thick,
            black,
            solid,
            postaction={
                decorate,
                decoration={
                    markings,
                    mark=at position 0.65 with {\arrow[scale=1.2]{stealth}}, 
                }
            }, 
            forget plot
        ] table[x index=0, y index=1] {data/pathData_dyn_example2.dat};
        \addlegendentry{Reference}

        \addplot[very thick, brown,
             mark=*, mark options={fill=brown}] table[x index=0, y index=1] {data/pathData_fesl_example2.dat};
        \addplot[
            very thick,
            brown,
            postaction={
                decorate,
                decoration={
                    markings,
                    mark=at position 0.4 with {\arrow[scale=1.5]{stealth}},
                    mark=at position 0.9 with {\arrow[scale=1.2]{stealth}}
                }
            }, 
            forget plot
        ] table[x index=0, y index=1] {data/pathData_fesl_example2.dat}; 
        \addlegendentry{F-ESL}

        \addplot[only marks, mark=*, blue] coordinates {(0.2743, 0.2500)};
        \node[anchor=north, blue] at (axis cs:0.2743, 0.2500) {$\mathbf{x}^{*}$};

	\end{axis}
\end{tikzpicture}

%% file: elements/system_drawing_example3.tex
\begin{tikzpicture}[every node/.style={draw,outer sep=0pt,thick}]

    \tikzset{
      ground/.style={
        draw=none,
        minimum width=0.6cm,
        minimum height=0.3cm,
        path picture={
          \fill[white]
            (path picture bounding box.south west) rectangle (path picture bounding box.north east);
          \fill[pattern=north east lines]
            (path picture bounding box.south west) rectangle (path picture bounding box.north east);
        }
      }
    }

    \coordinate (N1) at (-6, 0);
    \coordinate (N2) at (-3, 0);
    \coordinate (N3) at (0, 0);
    \coordinate (N4) at (3, 0);
    \coordinate (N5) at (6, 0);
    \coordinate (N6) at (-4, 1);
    \coordinate (N7) at (4, 1);
    \coordinate (N8) at (-2, 2);
    \coordinate (N9) at (2, 2);
    \coordinate (N10) at (0, 3);
    
    \def\barWidth{2mm} 

    \newcommand{\drawbar}[2]{%
        \path (#1); \pgfgetlastxy{\xa}{\ya}
        \path (#2); \pgfgetlastxy{\xb}{\yb}
        \pgfmathsetmacro{\dx}{\xb - \xa}
        \pgfmathsetmacro{\dy}{\yb - \ya}
        \pgfmathsetmacro{\angle}{atan2(\dy,\dx)}
        \pgfmathsetmacro{\len}{veclen(\dx,\dy)}
    
        \begin{scope}
            \path (#1);
            \coordinate (A) at ($(#1)+({\angle+90}:\barWidth/2)$);
            \coordinate (B) at ($(#2)+({\angle+90}:\barWidth/2)$);
            \coordinate (C) at ($(#2)+({\angle-90}:\barWidth/2)$);
            \coordinate (D) at ($(#1)+({\angle-90}:\barWidth/2)$);
            \draw[draw=black, fill=lightgray, very thick] (A) -- (B) -- (C) -- (D) -- cycle;
        \end{scope}
    }

    \drawbar{N1}{N6}
    \drawbar{N6}{N2}
    \drawbar{N4}{N7}
    \drawbar{N5}{N7}
    \drawbar{N6}{N8}
    \drawbar{N2}{N8}
    \drawbar{N3}{N8}
    \drawbar{N3}{N9}
    \drawbar{N4}{N9}
    \drawbar{N7}{N9}
    \drawbar{N8}{N10}
    \drawbar{N3}{N10}
    \drawbar{N9}{N10}
    
    \foreach \i in {1,2,3,4,5} {
        \draw [very thick] (N\i.north) -- + (0.3cm,-0.5cm) -- +(-0.3cm,-0.5cm) -- + (0cm,0cm);
    }
    
    \foreach \i in {1,2,3,4,5,6,7,8,9,10} {
        \node[anchor=center, draw, circle, fill=white, inner sep=3pt, very thick] at (N\i) {};
    }
    
    \foreach \i in {1,2,3,4,5} {
        \node (GRND\i) at ($(N\i)+(0,-0.5cm)$) [ground, anchor=north, draw=none] {};
        \draw[very thick] ($(GRND\i.north)-(4mm,0)$) -- ($(GRND\i.north)+(4mm,0)$);
    }

    \pgfmathsetmacro{\angle}{26.5651
}
    \pgfmathsetmacro{\sinang}{sin(\angle)}
    \pgfmathsetmacro{\cosang}{cos(\angle)}
    
    \draw [-latex, ultra thick] 
      ($(N6) + ({-1cm*\sinang}, {1cm*\cosang})$) -- (N6);
    \node [draw=none, anchor=west] at ($(N6)+(-1cm, 0.5cm)$) {$F_{w}$};
    
    \draw [-latex, ultra thick] 
      ($(N8) + ({-1cm*\sinang}, {1cm*\cosang})$) -- (N8);
    \node [draw=none, anchor=west] at ($(N8)+(-1cm, 0.5cm)$) {$F_{w}$};
    
    \draw [-latex, ultra thick] 
      ($(N10) + ({-1cm*\sinang}, {1cm*\cosang})$) -- (N10); 
    \node [draw=none, anchor=west] at ($(N10)+(-1cm, 0.5cm)$) {$F_{w}$};

    \draw [-latex, ultra thick] 
      (N10) -- ($(N10) + ({1cm*\sinang}, {1cm*\cosang})$);
    \node [draw=none, anchor=east] at ($(N10)+(1cm, 0.5cm)$) {$F_{w}$};
    
    \draw [-latex, ultra thick] 
      (N9) -- ($(N9) + ({1cm*\sinang}, {1cm*\cosang})$);
    \node [draw=none, anchor=east] at ($(N9)+(1cm, 0.5cm)$) {$F_{w}$};
    
    \draw [-latex, ultra thick] 
      (N7) -- ($(N7) + ({1cm*\sinang}, {1cm*\cosang})$); 
    \node [draw=none, anchor=east] at ($(N7)+(1cm, 0.5cm)$) {$F_{w}$};    


\end{tikzpicture}

%% file: elements/nodesData_example3.tex
\begin{tabular}{lrl}
\toprule
Node & $x$ & $y$ \\
\midrule
1  & $-6$ & $0$ \\
2  & $-3$ & $0$ \\
3  & $0$  & $0$ \\
4  & $3$  & $0$ \\
5  & $6$  & $0$ \\
6  & $-4$ & $1$ \\
7  & $4$  & $1$ \\
8  & $-2$ & $2$ \\
9  & $2$  & $2$ \\
10 & $0$  & $3$ \\
\bottomrule
\end{tabular}

%% file: elements/elementsData_example3.tex
\begin{tabular}{lll}
\toprule
Element & Node 1 & Node 2 \\
\midrule
1  & 1  & 6  \\
2  & 6  & 2  \\
3  & 4  & 7  \\
4  & 5  & 7  \\
5  & 6  & 8  \\
6  & 2  & 8  \\
7  & 3  & 8  \\
8  & 3  & 9  \\
9  & 4  & 9  \\
10 & 7  & 9  \\
11 & 8  & 10 \\
12 & 3  & 10 \\
13 & 9  & 10 \\
\bottomrule
\end{tabular}

%% file: elements/loadCases_def_table_example3.tex
\begin{tabular}{lccccccc}
\toprule
\multicolumn{8}{l}{
  \parbox{0.65\linewidth}{
    The loads are defined as \\[0.5ex]
    $p_i = \bar{p}_i \phi_i(t),$ where $\phi_i(t) = a_0 + \sum_{i=1}^3 a_i \sin(2\pi f_i t)$ \\[0.5ex]
    and $\bar{p}_1=\bar{p}_2=25 \times 10^5$ N.
  }
} \\
\midrule
$p_i(t)$ & $a_0$ & $a_1$ & $a_2$ & $a_3$ & $f_1$ & $f_2$ & $f_3$ \\
         &       &       &       &       & [Hz] & [Hz] & [Hz] \\
\midrule
$p_1(t)$ & 0.75 & 0.10 & 0.05 & 0.10 & 6.0 & 1.2 & 0.06 \\
$p_2(t)$ & 0.80 & 0.05 & 0.10 & 0.05 & 2.4 & 1.8 & 0.03 \\
\bottomrule
\end{tabular}

%% file: elements/loadsCases_example3_plot.tex


    

\begin{tikzpicture}
    \begin{axis}[
        xlabel={$t$ [s]},
        ylabel={$p_i(t)$ [N]},
        width=11cm,
        height=6cm,
        xmin=0, xmax=50, 
        xtick={0,10,20,30,40,50},
        grid
    ]

    \addplot[color=orange, thick] table[x index=0, y index=1] {data/LoadCase1Data_example3.dat};
    \addlegendentry{$p_1(t)$}    

    \addplot[color=blue, thick] table[x index=0, y index=1] {data/LoadCase2Data_example3.dat};
    \addlegendentry{$p_2(t)$}    
    
    \end{axis}
\end{tikzpicture}

%% file: elements/system_drawing_opt_example3.tex
\begin{tikzpicture}[every node/.style={draw,outer sep=0pt,thick}]

    \tikzset{
      ground/.style={
        draw=none,
        minimum width=0.6cm,
        minimum height=0.3cm,
        path picture={
          \fill[white]
            (path picture bounding box.south west) rectangle (path picture bounding box.north east);
          \fill[pattern=north east lines]
            (path picture bounding box.south west) rectangle (path picture bounding box.north east);
        }
      }
    }

    \coordinate (N1) at (-6, 0);
    \coordinate (N2) at (-3, 0);
    \coordinate (N3) at (0, 0);
    \coordinate (N4) at (3, 0);
    \coordinate (N5) at (6, 0);
    \coordinate (N6) at (-4, 1);
    \coordinate (N7) at (4, 1);
    \coordinate (N8) at (-2, 2);
    \coordinate (N9) at (2, 2);
    \coordinate (N10) at (0, 3);

    \newcommand{\drawbar}[3]{%
        \path (#1); \pgfgetlastxy{\xa}{\ya}
        \path (#2); \pgfgetlastxy{\xb}{\yb}
        \pgfmathsetmacro{\dx}{\xb - \xa}
        \pgfmathsetmacro{\dy}{\yb - \ya}
        \pgfmathsetmacro{\angle}{atan2(\dy,\dx)}
        \pgfmathsetmacro{\len}{veclen(\dx,\dy)}
    
        \begin{scope}
            \path (#1);
            \coordinate (A) at ($(#1)+({\angle+90}:#3/2)$);
            \coordinate (B) at ($(#2)+({\angle+90}:#3/2)$);
            \coordinate (C) at ($(#2)+({\angle-90}:#3/2)$);
            \coordinate (D) at ($(#1)+({\angle-90}:#3/2)$);
            \draw[draw=black, fill=lightgray, very thick] (A) -- (B) -- (C) -- (D) -- cycle;
        \end{scope}
    }

    \pgfmathsetmacro{\varA}{0.8375}
    \pgfmathsetmacro{\varB}{2.0000}
    \pgfmathsetmacro{\varC}{1.1830}
    \pgfmathsetmacro{\varD}{0.9241}
    \pgfmathsetmacro{\varE}{1.9385}
    \pgfmathsetmacro{\varF}{1.3851}
    \pgfmathsetmacro{\varG}{0.0155}
    
    \drawbar{N1}{N6} {\varA mm}
    \drawbar{N6}{N2} {\varB mm}
    \drawbar{N4}{N7} {\varB mm}
    \drawbar{N5}{N7} {\varA mm}
    \drawbar{N6}{N8} {\varC mm}
    \drawbar{N2}{N8} {\varD mm}
    \drawbar{N3}{N8} {\varE mm}
    \drawbar{N3}{N9} {\varE mm}
    \drawbar{N4}{N9} {\varD mm}
    \drawbar{N7}{N9} {\varC mm}
    \drawbar{N8}{N10}{\varF mm}
    \drawbar{N3}{N10}{\varG mm}
    \drawbar{N9}{N10}{\varF mm}
    
    \foreach \i in {1,2,3,4,5} {
        \draw [very thick] (N\i.north) -- + (0.3cm,-0.5cm) -- +(-0.3cm,-0.5cm) -- + (0cm,0cm);
    }
    
    \foreach \i in {1,2,3,4,5,6,7,8,9,10} {
        \node[anchor=center, draw, circle, fill=white, inner sep=3pt, very thick] at (N\i) {};
    }
    
    \foreach \i in {1,2,3,4,5} {
        \node (GRND\i) at ($(N\i)+(0,-0.5cm)$) [ground, anchor=north, draw=none] {};
        \draw[very thick] ($(GRND\i.north)-(4mm,0)$) -- ($(GRND\i.north)+(4mm,0)$);
    }
    

\end{tikzpicture}

%% file: elements/results_table_example3_1.tex
\begin{tabular}{llcccccccccc} \toprule
 &  & $N_{THA}$ & $x^{*}_1$ & $x^{*}_2$ & $x^{*}_3$ & $x^{*}_4$ & $x^{*}_5$ & $x^{*}_6$ & $x^{*}_7$ & $f_c(\mathbf{x}^{*})$ \\ 
\midrule
\eqref{eq:dynoptform_example3} & SQP & $35$ & $0.0541$ & $0.1292$ & $0.0764$ & $0.0597$ & $0.1253$ & $0.0895$ & $0.0010$ & $0.1783$ \\
\eqref{eq:statoptform_example3} & F-ESL & $4$  & $0.0541$ & $0.1292$ & $0.0764$ & $0.0597$ & $0.1253$ & $0.0895$ & $0.0010$ & $0.1783$ \\
\bottomrule
\end{tabular}

%% file: elements/results_table_example3_2.tex
\begin{tabular}{llcccccccc} \toprule
 &  
 & $g_{\max}$ 
 & $\frac{\partial \mathcal{L}(\mathbf{x}^*)}{\partial x_1}$ & $\frac{\partial \mathcal{L}(\mathbf{x}^*)}{\partial x_2}$ & $\frac{\partial \mathcal{L}(\mathbf{x}^*)}{\partial x_3}$ 
 & $\frac{\partial \mathcal{L}(\mathbf{x}^*)}{\partial x_4}$ & $\frac{\partial \mathcal{L}(\mathbf{x}^*)}{\partial x_5}$ & $\frac{\partial \mathcal{L}(\mathbf{x}^*)}{\partial x_6}$ & $\frac{\partial \mathcal{L}(\mathbf{x}^*)}{\partial x_7}$ \\
\midrule
\eqref{eq:dynoptform_example3} & SQP 
& $-3.6\mathrm{e}{-8}$ 
& $3.9\mathrm{e}{-9}$ & $3.8\mathrm{e}{-6}$ & $4.2\mathrm{e}{-7}$ 
& $-9.9\mathrm{e}{-8}$ & $3.3\mathrm{e}{-7}$ & $2.4\mathrm{e}{-7}$ & $5.9\mathrm{e}{-8}$ \\
\eqref{eq:statoptform_example3} & F-ESL 
& $-7.8\mathrm{e}{-9}$ 
& $2.0\mathrm{e}{-10}$ & $-1.0\mathrm{e}{-5}$ & $-2.3\mathrm{e}{-6}$ 
& $-4.8\mathrm{e}{-6}$ & $-9.3\mathrm{e}{-6}$ & $-2.3\mathrm{e}{-6}$ & $-1.3\mathrm{e}{-6}$ \\
\bottomrule
\end{tabular}

%% file: elements/sigma_plot_loadcase1_example3.tex
\begin{tikzpicture}
  \begin{groupplot}[
    group style={
      group size=2 by 3,
      horizontal sep=2.5cm,
      vertical sep=2.5cm,
    },
    width=6cm,
    height=4cm,
    xlabel={$t$ [s]},
    ylabel={$\sigma(t)$ [Pa]},
    xmin=30, xmax=50,
    ymin=-250e6, ymax=250e6,
    title style={yshift=2ex}
  ]

    \nextgroupplot[title={$x_1:$ bars 1 and 4}]
        \addplot[color=teal] table[x index=0, y index=1] {data/sig/sig1_example3.dat};
        \addlegendentry{$\sigma_1$}    
        \addplot[color=brown] table[x index=0, y index=1] {data/sig/sig4_example3.dat};
        \addlegendentry{$\sigma_4$}    
        \addplot[color=red, dashed, forget plot] table[x index=0, y index=1] {data/sig/sig_y_plus_example3.dat};
        \addplot[color=red, dashed, forget plot] table[x index=0, y index=1] {data/sig/sig_y_minus_example3.dat};
        \addplot[color=black, dashed, forget plot] table[x index=0, y index=1] {data/sig/sig1_buck_example3.dat};
    
    \nextgroupplot[title={$x_2:$ bars 2 and 3}] 
        \addplot[color=teal] table[x index=0, y index=1] {data/sig/sig2_example3.dat};
        \addlegendentry{$\sigma_2$}    
        \addplot[color=brown] table[x index=0, y index=1] {data/sig/sig3_example3.dat};
        \addlegendentry{$\sigma_3$}    
        \addplot[color=red, dashed, forget plot] table[x index=0, y index=1] {data/sig/sig_y_plus_example3.dat};
        \addplot[color=red, dashed, forget plot] table[x index=0, y index=1] {data/sig/sig_y_minus_example3.dat};
        \addplot[color=black, dashed, forget plot] table[x index=0, y index=1] {data/sig/sig2_buck_example3.dat};
    
    \nextgroupplot[title={$x_3:$ bars 5 and 10}] 
        \addplot[color=teal] table[x index=0, y index=1] {data/sig/sig5_example3.dat};
        \addlegendentry{$\sigma_5$}    
        \addplot[color=brown] table[x index=0, y index=1] {data/sig/sig10_example3.dat};
        \addlegendentry{$\sigma_{10}$}    
        \addplot[color=red, dashed, forget plot] table[x index=0, y index=1] {data/sig/sig_y_plus_example3.dat};
        \addplot[color=red, dashed, forget plot] table[x index=0, y index=1] {data/sig/sig_y_minus_example3.dat};
        \addplot[color=black, dashed, forget plot] table[x index=0, y index=1] {data/sig/sig5_buck_example3.dat};
    
    \nextgroupplot[title={$x_4:$ bars 6 and 9}] 
        \addplot[color=teal] table[x index=0, y index=1] {data/sig/sig6_example3.dat};
        \addlegendentry{$\sigma_6$}    
        \addplot[color=brown] table[x index=0, y index=1] {data/sig/sig9_example3.dat};
        \addlegendentry{$\sigma_9$}    
        \addplot[color=red, dashed, forget plot] table[x index=0, y index=1] {data/sig/sig_y_plus_example3.dat};
        \addplot[color=red, dashed, forget plot] table[x index=0, y index=1] {data/sig/sig_y_minus_example3.dat};
        \addplot[color=black, dashed, forget plot] table[x index=0, y index=1] {data/sig/sig6_buck_example3.dat};
    
    \nextgroupplot[title={$x_5:$ bars 7 and 8}] 
        \addplot[color=teal] table[x index=0, y index=1] {data/sig/sig7_example3.dat};
        \addlegendentry{$\sigma_7$}    
        \addplot[color=brown] table[x index=0, y index=1] {data/sig/sig8_example3.dat};
        \addlegendentry{$\sigma_8$}    
        \addplot[color=red, dashed, forget plot] table[x index=0, y index=1] {data/sig/sig_y_plus_example3.dat};
        \addplot[color=red, dashed, forget plot] table[x index=0, y index=1] {data/sig/sig_y_minus_example3.dat};
        \addplot[color=black, dashed, forget plot] table[x index=0, y index=1] {data/sig/sig7_buck_example3.dat};
    
    \nextgroupplot[title={$x_6:$ bars 11 and 13}] 
        \addplot[color=teal] table[x index=0, y index=1] {data/sig/sig11_example3.dat};
        \addlegendentry{$\sigma_{11}$}    
        \addplot[color=brown] table[x index=0, y index=1] {data/sig/sig13_example3.dat};
        \addlegendentry{$\sigma_{13}$}    
        \addplot[color=red, dashed, forget plot] table[x index=0, y index=1] {data/sig/sig_y_plus_example3.dat};
        \addplot[color=red, dashed, forget plot] table[x index=0, y index=1] {data/sig/sig_y_minus_example3.dat};
        \addplot[color=black, dashed, forget plot] table[x index=0, y index=1] {data/sig/sig11_buck_example3.dat};

  \end{groupplot}

  \path (current bounding box.south) ++(0,-0.5cm) node (legend) {
    \begin{tikzpicture}
    \begin{axis}[
            hide axis,
            xmin=0,
            xmax=1,
            ymin=0,
            ymax=1,
            legend columns=3,
            legend style={draw=none, column sep=1ex},
        ]

    \addlegendimage{red, thick, dashed}
    \addlegendentry{$\pm \sigma_y$}
    \addlegendimage{black, thick, dashed}
    \addlegendentry{$\sigma_{buck}$}
    \end{axis}
    \end{tikzpicture}
};

\end{tikzpicture}

%% file: elements/sigma_plot_loadcase2_example3.tex
\begin{tikzpicture}
  \begin{groupplot}[
    group style={
      group size=2 by 3,
      horizontal sep=2.5cm,
      vertical sep=2.5cm,
    },
    width=6cm,
    height=4cm,
    xlabel={$t$ [s]},
    ylabel={$\sigma(t)$ [Pa]},
    xmin=30, xmax=50,
    ymin=-250e6, ymax=250e6,
    title style={yshift=2ex}
  ]

    \nextgroupplot[title={$x_1:$ bars 1 and 4}] 
        \addplot[color=teal, thick] table[x index=0, y index=1] {data/sig2/sig1_2_example3.dat};
        \addlegendentry{$\sigma_1$}    
        \addplot[color=brown, thick] table[x index=0, y index=1] {data/sig2/sig4_2_example3.dat};
        \addlegendentry{$\sigma_4$}    
        \addplot[color=red, thick, dashed, forget plot] table[x index=0, y index=1] {data/sig2/sig_y_plus_2_example3.dat};
        \addplot[color=red, thick, dashed, forget plot] table[x index=0, y index=1] {data/sig2/sig_y_minus_2_example3.dat};
        \addplot[color=black, thick, dashed, forget plot] table[x index=0, y index=1] {data/sig2/sig1_buck_2_example3.dat};
    
    \nextgroupplot[title={$x_2:$ bars 2 and 3}] 
        \addplot[color=teal, thick] table[x index=0, y index=1] {data/sig2/sig2_2_example3.dat};
        \addlegendentry{$\sigma_2$}    
        \addplot[color=brown, thick] table[x index=0, y index=1] {data/sig2/sig3_2_example3.dat};
        \addlegendentry{$\sigma_3$}    
        \addplot[color=red, thick, dashed, forget plot] table[x index=0, y index=1] {data/sig2/sig_y_plus_2_example3.dat};
        \addplot[color=red, thick, dashed, forget plot] table[x index=0, y index=1] {data/sig2/sig_y_minus_2_example3.dat};
        \addplot[color=black, thick, dashed, forget plot] table[x index=0, y index=1] {data/sig2/sig2_buck_2_example3.dat};
    
    \nextgroupplot[title={$x_3:$ bars 5 and 10}] 
        \addplot[color=teal, thick] table[x index=0, y index=1] {data/sig2/sig5_2_example3.dat};
        \addlegendentry{$\sigma_5$}    
        \addplot[color=brown, thick] table[x index=0, y index=1] {data/sig2/sig10_2_example3.dat};
        \addlegendentry{$\sigma_{10}$}    
        \addplot[color=red, thick, dashed, forget plot] table[x index=0, y index=1] {data/sig2/sig_y_plus_2_example3.dat};
        \addplot[color=red, thick, dashed, forget plot] table[x index=0, y index=1] {data/sig2/sig_y_minus_2_example3.dat};
        \addplot[color=black, thick, dashed, forget plot] table[x index=0, y index=1] {data/sig2/sig5_buck_2_example3.dat};
    
    \nextgroupplot[title={$x_4:$ bars 6 and 9}] 
        \addplot[color=teal, thick] table[x index=0, y index=1] {data/sig2/sig6_2_example3.dat};
        \addlegendentry{$\sigma_6$}    
        \addplot[color=brown, thick] table[x index=0, y index=1] {data/sig2/sig9_2_example3.dat};
        \addlegendentry{$\sigma_9$}    
        \addplot[color=red, thick, dashed, forget plot] table[x index=0, y index=1] {data/sig2/sig_y_plus_2_example3.dat};
        \addplot[color=red, thick, dashed, forget plot] table[x index=0, y index=1] {data/sig2/sig_y_minus_2_example3.dat};
        \addplot[color=black, thick, dashed, forget plot] table[x index=0, y index=1] {data/sig2/sig6_buck_2_example3.dat};
    
    \nextgroupplot[title={$x_5:$ bars 7 and 8}] 
        \addplot[color=teal, thick] table[x index=0, y index=1] {data/sig2/sig7_2_example3.dat};
        \addlegendentry{$\sigma_7$}    
        \addplot[color=brown, thick] table[x index=0, y index=1] {data/sig2/sig8_2_example3.dat};
        \addlegendentry{$\sigma_8$}    
        \addplot[color=red, thick, dashed, forget plot] table[x index=0, y index=1] {data/sig2/sig_y_plus_2_example3.dat};
        \addplot[color=red, thick, dashed, forget plot] table[x index=0, y index=1] {data/sig2/sig_y_minus_2_example3.dat};
        \addplot[color=black, thick, dashed, forget plot] table[x index=0, y index=1] {data/sig2/sig7_buck_2_example3.dat};
    
    \nextgroupplot[title={$x_6:$ bars 11 and 13}] 
        \addplot[color=teal, thick] table[x index=0, y index=1] {data/sig2/sig11_2_example3.dat};
        \addlegendentry{$\sigma_{11}$}    
        \addplot[color=brown, thick] table[x index=0, y index=1] {data/sig2/sig13_2_example3.dat};
        \addlegendentry{$\sigma_{13}$}    
        \addplot[color=red, thick, dashed, forget plot] table[x index=0, y index=1] {data/sig2/sig_y_plus_2_example3.dat};
        \addplot[color=red, thick, dashed, forget plot] table[x index=0, y index=1] {data/sig2/sig_y_minus_2_example3.dat};
        \addplot[color=black, thick, dashed, forget plot] table[x index=0, y index=1] {data/sig2/sig11_buck_2_example3.dat};


  \end{groupplot}

  \path (current bounding box.south) ++(0,-0.5cm) node (legend) {
    \begin{tikzpicture}
    \begin{axis}[
            hide axis,
            xmin=0,
            xmax=1,
            ymin=0,
            ymax=1,
            legend columns=3,
            legend style={draw=none, column sep=1ex},
        ]

    \addlegendimage{red, thick, dashed}
    \addlegendentry{$\pm \sigma_y$}
    \addlegendimage{black, thick, dashed}
    \addlegendentry{$\sigma_{buck}$}
    \end{axis}
    \end{tikzpicture}
};

\end{tikzpicture}

%% file: manuscript.bbl
\begin{thebibliography}{55}
\expandafter\ifx\csname natexlab\endcsname\relax\def\natexlab#1{#1}\fi
\providecommand{\url}[1]{\texttt{#1}}
\providecommand{\href}[2]{#2}
\providecommand{\path}[1]{#1}
\providecommand{\DOIprefix}{doi:}
\providecommand{\ArXivprefix}{arXiv:}
\providecommand{\URLprefix}{URL: }
\providecommand{\Pubmedprefix}{pmid:}
\providecommand{\doi}[1]{\href{http://dx.doi.org/#1}{\path{#1}}}
\providecommand{\Pubmed}[1]{\href{pmid:#1}{\path{#1}}}
\providecommand{\bibinfo}[2]{#2}
\ifx\xfnm\relax \def\xfnm[#1]{\unskip,\space#1}\fi
\bibitem[{Amir et~al.(2008)Amir, Kirsch \& Sheinman}]{amir2008efficient}
\bibinfo{author}{Amir, O.}, \bibinfo{author}{Kirsch, U.}, \& \bibinfo{author}{Sheinman, I.} (\bibinfo{year}{2008}).
\newblock \bibinfo{title}{Efficient non-linear reanalysis of skeletal structures using combined approximations}.
\newblock {\it \bibinfo{journal}{International Journal for Numerical Methods in Engineering}\/},  {\it \bibinfo{volume}{73}\/}, \bibinfo{pages}{1328--1346}.
\bibitem[{Behrou \& Guest(2017)}]{behrou2017topology}
\bibinfo{author}{Behrou, R.}, \& \bibinfo{author}{Guest, J.~K.} (\bibinfo{year}{2017}).
\newblock \bibinfo{title}{Topology optimization for transient response of structures subjected to dynamic loads}.
\newblock In {\it \bibinfo{booktitle}{18th AIAA/ISSMO multidisciplinary analysis and optimization conference}\/} (p. \bibinfo{pages}{3657}).
\bibitem[{Blasques \& Stolpe(2012)}]{blasques2012multi}
\bibinfo{author}{Blasques, J.~P.}, \& \bibinfo{author}{Stolpe, M.} (\bibinfo{year}{2012}).
\newblock \bibinfo{title}{Multi-material topology optimization of laminated composite beam cross sections}.
\newblock {\it \bibinfo{journal}{Composite Structures}\/},  {\it \bibinfo{volume}{94}\/}, \bibinfo{pages}{3278--3289}.
\bibitem[{Choi et~al.(2018)Choi, Lee, Yoon, Han \& Park}]{choi2018structural}
\bibinfo{author}{Choi, W.-H.}, \bibinfo{author}{Lee, Y.}, \bibinfo{author}{Yoon, J.-M.}, \bibinfo{author}{Han, Y.-H.}, \& \bibinfo{author}{Park, G.-J.} (\bibinfo{year}{2018}).
\newblock \bibinfo{title}{Structural optimization for roof crush test using an enforced displacement method}.
\newblock {\it \bibinfo{journal}{International Journal of Automotive Technology}\/},  {\it \bibinfo{volume}{19}\/}, \bibinfo{pages}{291--299}.
\bibitem[{Choi \& Park(1999)}]{choi1999transformation}
\bibinfo{author}{Choi, W.-S.}, \& \bibinfo{author}{Park, G.} (\bibinfo{year}{1999}).
\newblock \bibinfo{title}{Transformation of dynamic loads into equivalent static loads based on modal analysis}.
\newblock {\it \bibinfo{journal}{International Journal for Numerical Methods in Engineering}\/},  {\it \bibinfo{volume}{46}\/}, \bibinfo{pages}{29--43}.
\bibitem[{Choi \& Park(2002)}]{choi2002structural}
\bibinfo{author}{Choi, W.-S.}, \& \bibinfo{author}{Park, G.-J.} (\bibinfo{year}{2002}).
\newblock \bibinfo{title}{Structural optimization using equivalent static loads at all time intervals}.
\newblock {\it \bibinfo{journal}{Computer Methods in Applied Mechanics and Engineering}\/},  {\it \bibinfo{volume}{191}\/}, \bibinfo{pages}{2105--2122}.
\bibitem[{Chopra(2007)}]{chopra2007dynamics}
\bibinfo{author}{Chopra, A.} (\bibinfo{year}{2007}).
\newblock {\it \bibinfo{title}{Dynamics of Structures}\/}.
\newblock Prentice-Hall international series in civil engineering and engineering mechanics.
\newblock \bibinfo{publisher}{Pearson Education}.
\newblock \URLprefix \url{https://books.google.co.il/books?id=0dU1bDaRyP4C}.
\bibitem[{Comi et~al.(2016)Comi, Corigliano, Zega \& Zerbini}]{comi2016non}
\bibinfo{author}{Comi, C.}, \bibinfo{author}{Corigliano, A.}, \bibinfo{author}{Zega, V.}, \& \bibinfo{author}{Zerbini, S.} (\bibinfo{year}{2016}).
\newblock \bibinfo{title}{Non linear response and optimization of a new z-axis resonant micro-accelerometer}.
\newblock {\it \bibinfo{journal}{Mechatronics}\/},  {\it \bibinfo{volume}{40}\/}, \bibinfo{pages}{235--243}.
\bibitem[{D{\'\i}aaz \& Kikuchi(1992)}]{diaaz1992solutions}
\bibinfo{author}{D{\'\i}aaz, A.~R.}, \& \bibinfo{author}{Kikuchi, N.} (\bibinfo{year}{1992}).
\newblock \bibinfo{title}{Solutions to shape and topology eigenvalue optimization problems using a homogenization method}.
\newblock {\it \bibinfo{journal}{International Journal for Numerical Methods in Engineering}\/},  {\it \bibinfo{volume}{35}\/}, \bibinfo{pages}{1487--1502}.
\bibitem[{Ferrari et~al.(2018)Ferrari, Lazarov \& Sigmund}]{ferrari2018eigenvalue}
\bibinfo{author}{Ferrari, F.}, \bibinfo{author}{Lazarov, B.~S.}, \& \bibinfo{author}{Sigmund, O.} (\bibinfo{year}{2018}).
\newblock \bibinfo{title}{Eigenvalue topology optimization via efficient multilevel solution of the frequency response}.
\newblock {\it \bibinfo{journal}{International Journal for Numerical Methods in Engineering}\/},  {\it \bibinfo{volume}{115}\/}, \bibinfo{pages}{872--892}.
\bibitem[{Giannini et~al.(2022)Giannini, Aage \& Braghin}]{giannini2022topology}
\bibinfo{author}{Giannini, D.}, \bibinfo{author}{Aage, N.}, \& \bibinfo{author}{Braghin, F.} (\bibinfo{year}{2022}).
\newblock \bibinfo{title}{Topology optimization of mems resonators with target eigenfrequencies and modes}.
\newblock {\it \bibinfo{journal}{European Journal of Mechanics-A/Solids}\/},  {\it \bibinfo{volume}{91}\/}, \bibinfo{pages}{104352}.
\bibitem[{Giannini et~al.(2020)Giannini, Braghin \& Aage}]{giannini2020topology}
\bibinfo{author}{Giannini, D.}, \bibinfo{author}{Braghin, F.}, \& \bibinfo{author}{Aage, N.} (\bibinfo{year}{2020}).
\newblock \bibinfo{title}{Topology optimization of 2d in-plane single mass mems gyroscopes}.
\newblock {\it \bibinfo{journal}{Structural and Multidisciplinary Optimization}\/},  {\it \bibinfo{volume}{62}\/}, \bibinfo{pages}{2069--2089}.
\bibitem[{Hermansen \& Lund(2024)}]{hermansen2024multi}
\bibinfo{author}{Hermansen, S.~M.}, \& \bibinfo{author}{Lund, E.} (\bibinfo{year}{2024}).
\newblock \bibinfo{title}{Multi-material and thickness optimization of a wind turbine blade root section}.
\newblock {\it \bibinfo{journal}{Structural and Multidisciplinary Optimization}\/},  {\it \bibinfo{volume}{67}\/}, \bibinfo{pages}{107}.
\bibitem[{Hudson(1995)}]{hudson1995dynamics}
\bibinfo{author}{Hudson, D.~E.} (\bibinfo{year}{1995}).
\newblock \bibinfo{title}{Dynamics of structures: Theory and applications to earthquake engineering, by anil k. chopra, prentice-hall, englewood cliffs, nj, 1995. no. of pages: xxviii+ 761, isbn 0-13-855214-2}.
\bibitem[{Idels \& Lavan(2020)}]{idels2020performance}
\bibinfo{author}{Idels, O.}, \& \bibinfo{author}{Lavan, O.} (\bibinfo{year}{2020}).
\newblock \bibinfo{title}{Performance based formal optimized seismic design of steel moment resisting frames}.
\newblock {\it \bibinfo{journal}{Computers \& Structures}\/},  {\it \bibinfo{volume}{235}\/}, \bibinfo{pages}{106269}.
\bibitem[{Idels \& Lavan(2023)}]{idels2023optimization}
\bibinfo{author}{Idels, O.}, \& \bibinfo{author}{Lavan, O.} (\bibinfo{year}{2023}).
\newblock \bibinfo{title}{Optimization-based seismic design of irregular self-centering moment resisting frames with ed bars or fluid viscous dampers}.
\newblock {\it \bibinfo{journal}{Structural and Multidisciplinary Optimization}\/},  {\it \bibinfo{volume}{66}\/}, \bibinfo{pages}{192}.
\bibitem[{James \& Waisman(2015)}]{james2015topology}
\bibinfo{author}{James, K.~A.}, \& \bibinfo{author}{Waisman, H.} (\bibinfo{year}{2015}).
\newblock \bibinfo{title}{Topology optimization of viscoelastic structures using a time-dependent adjoint method}.
\newblock {\it \bibinfo{journal}{Computer Methods in Applied Mechanics and Engineering}\/},  {\it \bibinfo{volume}{285}\/}, \bibinfo{pages}{166--187}.
\bibitem[{Jeong et~al.(2008)Jeong, Yi, Kan, Nagabhushana \& Park}]{jeong2008structural}
\bibinfo{author}{Jeong, S.}, \bibinfo{author}{Yi, S.}, \bibinfo{author}{Kan, C.}, \bibinfo{author}{Nagabhushana, V.}, \& \bibinfo{author}{Park, G.} (\bibinfo{year}{2008}).
\newblock \bibinfo{title}{Structural optimization of an automobile roof structure using equivalent static loads}.
\newblock {\it \bibinfo{journal}{Proceedings of the Institution of Mechanical Engineers, Part D: Journal of Automobile Engineering}\/},  {\it \bibinfo{volume}{222}\/}, \bibinfo{pages}{1985--1995}.
\bibitem[{Jeong et~al.(2010)Jeong, Yoon, Xu \& Park}]{jeong2010non}
\bibinfo{author}{Jeong, S.}, \bibinfo{author}{Yoon, S.}, \bibinfo{author}{Xu, S.}, \& \bibinfo{author}{Park, G.} (\bibinfo{year}{2010}).
\newblock \bibinfo{title}{Non-linear dynamic response structural optimization of an automobile frontal structure using equivalent static loads}.
\newblock {\it \bibinfo{journal}{Proceedings of the Institution of Mechanical Engineers, Part D: Journal of Automobile Engineering}\/},  {\it \bibinfo{volume}{224}\/}, \bibinfo{pages}{489--501}.
\bibitem[{Kang et~al.(2001)Kang, Choi \& Park}]{kang2001structural}
\bibinfo{author}{Kang, B.}, \bibinfo{author}{Choi, W.}, \& \bibinfo{author}{Park, G.} (\bibinfo{year}{2001}).
\newblock \bibinfo{title}{Structural optimization under equivalent static loads transformed from dynamic loads based on displacement}.
\newblock {\it \bibinfo{journal}{Computers \& Structures}\/},  {\it \bibinfo{volume}{79}\/}, \bibinfo{pages}{145--154}.
\bibitem[{Karev et~al.(2017)Karev, Harzheim, Immel \& Erzgr{\"a}ber}]{karev2017comparison}
\bibinfo{author}{Karev, A.}, \bibinfo{author}{Harzheim, L.}, \bibinfo{author}{Immel, R.}, \& \bibinfo{author}{Erzgr{\"a}ber, M.} (\bibinfo{year}{2017}).
\newblock \bibinfo{title}{{Comparison of different formulations of a front hood free sizing optimization problem using the ESL-method}}.
\newblock In {\it \bibinfo{booktitle}{World Congress of Structural and Multidisciplinary Optimisation}\/} (pp. \bibinfo{pages}{933--951}).
\newblock \bibinfo{organization}{Springer}.
\bibitem[{Karev et~al.(2019)Karev, Harzheim, Immel \& Erzgr{\"a}ber}]{karev2019free}
\bibinfo{author}{Karev, A.}, \bibinfo{author}{Harzheim, L.}, \bibinfo{author}{Immel, R.}, \& \bibinfo{author}{Erzgr{\"a}ber, M.} (\bibinfo{year}{2019}).
\newblock \bibinfo{title}{{Free sizing optimization of a front hood using the ESL method: overcoming challenges and traps}}.
\newblock {\it \bibinfo{journal}{Structural and Multidisciplinary Optimization}\/},  {\it \bibinfo{volume}{60}\/}, \bibinfo{pages}{1687--1707}.
\bibitem[{Kennedy \& Martins(2014)}]{kennedy2014parallel}
\bibinfo{author}{Kennedy, G.~J.}, \& \bibinfo{author}{Martins, J.~R.} (\bibinfo{year}{2014}).
\newblock \bibinfo{title}{A parallel aerostructural optimization framework for aircraft design studies}.
\newblock {\it \bibinfo{journal}{Structural and Multidisciplinary Optimization}\/},  {\it \bibinfo{volume}{50}\/}, \bibinfo{pages}{1079--1101}.
\bibitem[{Kirsch(1991)}]{kirsch1991reduced}
\bibinfo{author}{Kirsch, U.} (\bibinfo{year}{1991}).
\newblock \bibinfo{title}{Reduced basis approximations of structural displacements for optimal design}.
\newblock {\it \bibinfo{journal}{AIAA journal}\/},  {\it \bibinfo{volume}{29}\/}, \bibinfo{pages}{1751--1758}.
\bibitem[{Kirsch(2010)}]{kirsch2010reanalysis}
\bibinfo{author}{Kirsch, U.} (\bibinfo{year}{2010}).
\newblock \bibinfo{title}{Reanalysis and sensitivity reanalysis by combined approximations}.
\newblock {\it \bibinfo{journal}{Structural and Multidisciplinary Optimization}\/},  {\it \bibinfo{volume}{40}\/}, \bibinfo{pages}{1--15}.
\bibitem[{Kirsch \& Bogomolni(2004)}]{kirsch2004procedures}
\bibinfo{author}{Kirsch, U.}, \& \bibinfo{author}{Bogomolni, M.} (\bibinfo{year}{2004}).
\newblock \bibinfo{title}{Procedures for approximate eigenproblem reanalysis of structures}.
\newblock {\it \bibinfo{journal}{International Journal for Numerical Methods in Engineering}\/},  {\it \bibinfo{volume}{60}\/}, \bibinfo{pages}{1969--1986}.
\bibitem[{Kirsch et~al.(2006)Kirsch, Bogomolni \& Sheinman}]{kirsch2006nonlinear}
\bibinfo{author}{Kirsch, U.}, \bibinfo{author}{Bogomolni, M.}, \& \bibinfo{author}{Sheinman, I.} (\bibinfo{year}{2006}).
\newblock \bibinfo{title}{Nonlinear dynamic reanalysis of structures by combined approximations}.
\newblock {\it \bibinfo{journal}{Computer Methods in Applied Mechanics and Engineering}\/},  {\it \bibinfo{volume}{195}\/}, \bibinfo{pages}{4420--4432}.
\bibitem[{Kirsch et~al.(2007)Kirsch, Bogomolni \& Sheinman}]{kirsch2007efficient}
\bibinfo{author}{Kirsch, U.}, \bibinfo{author}{Bogomolni, M.}, \& \bibinfo{author}{Sheinman, I.} (\bibinfo{year}{2007}).
\newblock \bibinfo{title}{Efficient dynamic reanalysis of structures}.
\newblock {\it \bibinfo{journal}{Journal of Structural Engineering}\/},  {\it \bibinfo{volume}{133}\/}, \bibinfo{pages}{440--448}.
\bibitem[{Lee \& Park(2015)}]{lee2015nonlinear}
\bibinfo{author}{Lee, H.-A.}, \& \bibinfo{author}{Park, G.~J.} (\bibinfo{year}{2015}).
\newblock \bibinfo{title}{Nonlinear dynamic response topology optimization using the equivalent static loads method}.
\newblock {\it \bibinfo{journal}{Computer Methods in Applied Mechanics and Engineering}\/},  {\it \bibinfo{volume}{283}\/}, \bibinfo{pages}{956--970}.
\bibitem[{Li et~al.(2021)Li, Sigmund, Jensen \& Aage}]{li2021reduced}
\bibinfo{author}{Li, Q.}, \bibinfo{author}{Sigmund, O.}, \bibinfo{author}{Jensen, J.~S.}, \& \bibinfo{author}{Aage, N.} (\bibinfo{year}{2021}).
\newblock \bibinfo{title}{Reduced-order methods for dynamic problems in topology optimization: A comparative study}.
\newblock {\it \bibinfo{journal}{Computer Methods in Applied Mechanics and Engineering}\/},  {\it \bibinfo{volume}{387}\/}, \bibinfo{pages}{114149}.
\bibitem[{Li et~al.(2024)Li, Yin, Guo, Li \& Wang}]{li2023efficient}
\bibinfo{author}{Li, S.}, \bibinfo{author}{Yin, J.}, \bibinfo{author}{Guo, D.}, \bibinfo{author}{Li, G.}, \& \bibinfo{author}{Wang, H.} (\bibinfo{year}{2024}).
\newblock \bibinfo{title}{An efficient online successive reanalysis method for dynamic topology optimization}.
\newblock {\it \bibinfo{journal}{Advances in Engineering Software}\/},  {\it \bibinfo{volume}{191}\/}, \bibinfo{pages}{103625}.
\bibitem[{Long et~al.(2021)Long, Yang, Saeed, Tian, Wen \& Wang}]{long2021topology}
\bibinfo{author}{Long, K.}, \bibinfo{author}{Yang, X.}, \bibinfo{author}{Saeed, N.}, \bibinfo{author}{Tian, R.}, \bibinfo{author}{Wen, P.}, \& \bibinfo{author}{Wang, X.} (\bibinfo{year}{2021}).
\newblock \bibinfo{title}{Topology optimization of transient problem with maximum dynamic response constraint using soar scheme}.
\newblock {\it \bibinfo{journal}{Frontiers of Mechanical Engineering}\/},  {\it \bibinfo{volume}{16}\/}, \bibinfo{pages}{593--606}.
\bibitem[{Martins et~al.(2019)Martins, Sim{\~o}es \& Negr{\~a}o}]{martins2019optimization}
\bibinfo{author}{Martins, A.~M.}, \bibinfo{author}{Sim{\~o}es, L.~M.}, \& \bibinfo{author}{Negr{\~a}o, J.~H.} (\bibinfo{year}{2019}).
\newblock \bibinfo{title}{Optimization of concrete cable-stayed bridges under seismic action}.
\newblock {\it \bibinfo{journal}{Computers \& Structures}\/},  {\it \bibinfo{volume}{222}\/}, \bibinfo{pages}{36--47}.
\bibitem[{Marzok \& Lavan(2022)}]{marzok2022topology}
\bibinfo{author}{Marzok, A.}, \& \bibinfo{author}{Lavan, O.} (\bibinfo{year}{2022}).
\newblock \bibinfo{title}{Topology optimization of multiple-rocking concentrically braced frames subjected to earthquakes}.
\newblock {\it \bibinfo{journal}{Structural and Multidisciplinary Optimization}\/},  {\it \bibinfo{volume}{65}\/}, \bibinfo{pages}{104}.
\bibitem[{Niu et~al.(2018)Niu, He, Shan \& Yang}]{niu2018objective}
\bibinfo{author}{Niu, B.}, \bibinfo{author}{He, X.}, \bibinfo{author}{Shan, Y.}, \& \bibinfo{author}{Yang, R.} (\bibinfo{year}{2018}).
\newblock \bibinfo{title}{On objective functions of minimizing the vibration response of continuum structures subjected to external harmonic excitation}.
\newblock {\it \bibinfo{journal}{Structural and Multidisciplinary Optimization}\/},  {\it \bibinfo{volume}{57}\/}, \bibinfo{pages}{2291--2307}.
\bibitem[{Olhoff \& Du(2016)}]{olhoff2016generalized}
\bibinfo{author}{Olhoff, N.}, \& \bibinfo{author}{Du, J.} (\bibinfo{year}{2016}).
\newblock \bibinfo{title}{Generalized incremental frequency method for topological designof continuum structures for minimum dynamic compliance subject to forced vibration at a prescribed low or high value of the excitation frequency}.
\newblock {\it \bibinfo{journal}{Structural and Multidisciplinary Optimization}\/},  {\it \bibinfo{volume}{54}\/}, \bibinfo{pages}{1113--1141}.
\bibitem[{Park \& Kang(2003)}]{park2003validation}
\bibinfo{author}{Park, G.~J.}, \& \bibinfo{author}{Kang, B.} (\bibinfo{year}{2003}).
\newblock \bibinfo{title}{Validation of a structural optimization algorithm transforming dynamic loads into equivalent static loads}.
\newblock {\it \bibinfo{journal}{Journal of Optimization Theory and Applications}\/},  {\it \bibinfo{volume}{118}\/}, \bibinfo{pages}{191--200}.
\bibitem[{Park et~al.(2005)Park, Lee \& Park}]{park2005structural}
\bibinfo{author}{Park, K.~J.}, \bibinfo{author}{Lee, J.~N.}, \& \bibinfo{author}{Park, G.~J.} (\bibinfo{year}{2005}).
\newblock \bibinfo{title}{Structural shape optimization using equivalent static loads transformed from dynamic loads}.
\newblock {\it \bibinfo{journal}{International Journal for Numerical Methods in Engineering}\/},  {\it \bibinfo{volume}{63}\/}, \bibinfo{pages}{589--602}.
\bibitem[{Pavese et~al.(2017)Pavese, Tibaldi, Zahle \& Kim}]{pavese2017aeroelastic}
\bibinfo{author}{Pavese, C.}, \bibinfo{author}{Tibaldi, C.}, \bibinfo{author}{Zahle, F.}, \& \bibinfo{author}{Kim, T.} (\bibinfo{year}{2017}).
\newblock \bibinfo{title}{Aeroelastic multidisciplinary design optimization of a swept wind turbine blade}.
\newblock {\it \bibinfo{journal}{Wind Energy}\/},  {\it \bibinfo{volume}{20}\/}, \bibinfo{pages}{1941--1953}.
\bibitem[{Pedersen(2004)}]{pedersen2004crashworthiness}
\bibinfo{author}{Pedersen, C.~B.} (\bibinfo{year}{2004}).
\newblock \bibinfo{title}{Crashworthiness design of transient frame structures using topology optimization}.
\newblock {\it \bibinfo{journal}{Computer Methods in Applied Mechanics and Engineering}\/},  {\it \bibinfo{volume}{193}\/}, \bibinfo{pages}{653--678}.
\bibitem[{Pedersen(2000)}]{pedersen2000maximization}
\bibinfo{author}{Pedersen, N.~L.} (\bibinfo{year}{2000}).
\newblock \bibinfo{title}{Maximization of eigenvalues using topology optimization}.
\newblock {\it \bibinfo{journal}{Structural and Multidisciplinary Optimization}\/},  {\it \bibinfo{volume}{20}\/}, \bibinfo{pages}{2--11}.
\bibitem[{Przemieniecki(1985)}]{przemieniecki1985theory}
\bibinfo{author}{Przemieniecki, J.~S.} (\bibinfo{year}{1985}).
\newblock {\it \bibinfo{title}{Theory of matrix structural analysis}\/}.
\newblock \bibinfo{publisher}{Courier Corporation}.
\bibitem[{Somerville(1997)}]{somerville1997development}
\bibinfo{author}{Somerville, P.~G.} (\bibinfo{year}{1997}).
\newblock {\it \bibinfo{title}{Development of ground motion time histories for phase 2 of the FEMA/SAC steel project}\/}.
\newblock \bibinfo{type}{Technical Report} SAC Joint Venture.
\bibitem[{Stolpe(2014)}]{stolpe2014equivalent}
\bibinfo{author}{Stolpe, M.} (\bibinfo{year}{2014}).
\newblock \bibinfo{title}{On the equivalent static loads approach for dynamic response structural optimization}.
\newblock {\it \bibinfo{journal}{Structural and Multidisciplinary Optimization}\/},  {\it \bibinfo{volume}{50}\/}, \bibinfo{pages}{921--926}.
\bibitem[{Stolpe(2023)}]{stolpe2023diesl}
\bibinfo{author}{Stolpe, M.} (\bibinfo{year}{2023}).
\newblock \bibinfo{title}{When is diesl= esl for linear dynamic response structural optimization?}
\newblock {\it \bibinfo{journal}{Structural and Multidisciplinary Optimization}\/},  {\it \bibinfo{volume}{66}\/}, \bibinfo{pages}{73}.
\bibitem[{Stolpe \& Pollini(2023)}]{stolpe2023first}
\bibinfo{author}{Stolpe, M.}, \& \bibinfo{author}{Pollini, N.} (\bibinfo{year}{2023}).
\newblock \bibinfo{title}{A first-order equivalent static loads algorithm for optimization of nonlinear static response}.
\newblock {\it \bibinfo{journal}{Advances in Engineering Software}\/},  {\it \bibinfo{volume}{182}\/}, \bibinfo{pages}{103462}.
\bibitem[{Stolpe et~al.(2018)Stolpe, Verbart \& Rojas-Labanda}]{stolpe2018equivalent}
\bibinfo{author}{Stolpe, M.}, \bibinfo{author}{Verbart, A.}, \& \bibinfo{author}{Rojas-Labanda, S.} (\bibinfo{year}{2018}).
\newblock \bibinfo{title}{The equivalent static loads method for structural optimization does not in general generate optimal designs}.
\newblock {\it \bibinfo{journal}{Structural and Multidisciplinary Optimization}\/},  {\it \bibinfo{volume}{58}\/}, \bibinfo{pages}{139--154}.
\bibitem[{Torstenfelt \& Klarbring(2006)}]{torstenfelt2006structural}
\bibinfo{author}{Torstenfelt, B.}, \& \bibinfo{author}{Klarbring, A.} (\bibinfo{year}{2006}).
\newblock \bibinfo{title}{Structural optimization of modular product families with application to car space frame structures}.
\newblock {\it \bibinfo{journal}{Structural and Multidisciplinary optimization}\/},  {\it \bibinfo{volume}{32}\/}, \bibinfo{pages}{133--140}.
\bibitem[{Triller et~al.(2022)Triller, Immel \& Harzheim}]{triller2022difference}
\bibinfo{author}{Triller, J.}, \bibinfo{author}{Immel, R.}, \& \bibinfo{author}{Harzheim, L.} (\bibinfo{year}{2022}).
\newblock \bibinfo{title}{Difference-based equivalent static load method with adaptive time selection and local stiffness adaption}.
\newblock {\it \bibinfo{journal}{Structural and Multidisciplinary Optimization}\/},  {\it \bibinfo{volume}{65}\/}, \bibinfo{pages}{89}.
\bibitem[{Triller et~al.(2021)Triller, Immel, Timmer \& Harzheim}]{Triller2021}
\bibinfo{author}{Triller, J.}, \bibinfo{author}{Immel, R.}, \bibinfo{author}{Timmer, A.}, \& \bibinfo{author}{Harzheim, L.} (\bibinfo{year}{2021}).
\newblock \bibinfo{title}{{The difference-based equivalent static load method: an improvement of the ESL method’s nonlinear approximation quality}}.
\newblock {\it \bibinfo{journal}{Structural and Multidisciplinary Optimization}\/},  {\it \bibinfo{volume}{63}\/}, \bibinfo{pages}{2705--2720}.
\bibitem[{Tromme et~al.(2016)Tromme, Sonneville, Br{\"u}ls \& Duysinx}]{tromme2016equivalent}
\bibinfo{author}{Tromme, E.}, \bibinfo{author}{Sonneville, V.}, \bibinfo{author}{Br{\"u}ls, O.}, \& \bibinfo{author}{Duysinx, P.} (\bibinfo{year}{2016}).
\newblock \bibinfo{title}{On the equivalent static load method for flexible multibody systems described with a nonlinear finite element formalism}.
\newblock {\it \bibinfo{journal}{International Journal for Numerical Methods in Engineering}\/},  {\it \bibinfo{volume}{108}\/}, \bibinfo{pages}{646--664}.
\bibitem[{Tromme et~al.(2018)Tromme, Sonneville, Guest \& Br{\"u}ls}]{tromme2018system}
\bibinfo{author}{Tromme, E.}, \bibinfo{author}{Sonneville, V.}, \bibinfo{author}{Guest, J.~K.}, \& \bibinfo{author}{Br{\"u}ls, O.} (\bibinfo{year}{2018}).
\newblock \bibinfo{title}{System-wise equivalent static loads for the design of flexible mechanisms}.
\newblock {\it \bibinfo{journal}{Computer Methods in Applied Mechanics and Engineering}\/},  {\it \bibinfo{volume}{329}\/}, \bibinfo{pages}{312--331}.
\bibitem[{Verbart \& Stolpe(2018)}]{verbart2018working}
\bibinfo{author}{Verbart, A.}, \& \bibinfo{author}{Stolpe, M.} (\bibinfo{year}{2018}).
\newblock \bibinfo{title}{A working-set approach for sizing optimization of frame-structures subjected to time-dependent constraints}.
\newblock {\it \bibinfo{journal}{Structural and Multidisciplinary Optimization}\/},  {\it \bibinfo{volume}{58}\/}, \bibinfo{pages}{1367--1382}.
\bibitem[{Zhao \& Wang(2016)}]{zhao2016dynamic}
\bibinfo{author}{Zhao, J.}, \& \bibinfo{author}{Wang, C.} (\bibinfo{year}{2016}).
\newblock \bibinfo{title}{Dynamic response topology optimization in the time domain using model reduction method}.
\newblock {\it \bibinfo{journal}{Structural and Multidisciplinary Optimization}\/},  {\it \bibinfo{volume}{53}\/}, \bibinfo{pages}{101--114}.
\bibitem[{Zordan et~al.(2010)Zordan, Briseghella \& Mazzarolo}]{zordan2010bridge}
\bibinfo{author}{Zordan, T.}, \bibinfo{author}{Briseghella, B.}, \& \bibinfo{author}{Mazzarolo, E.} (\bibinfo{year}{2010}).
\newblock \bibinfo{title}{Bridge structural optimization through step-by-step evolutionary process}.
\newblock {\it \bibinfo{journal}{Structural Engineering International}\/},  {\it \bibinfo{volume}{20}\/}, \bibinfo{pages}{72--78}.

\end{thebibliography}
